\documentclass[11pt]{article}

\usepackage{graphics,enumitem,epsfig,textcomp}
\usepackage{amsfonts,amsmath,amssymb,amsthm}
\usepackage{euscript,color,mathrsfs}

\usepackage{url}
\usepackage{hyperref}
\hypersetup{
    colorlinks=true,
    linkcolor=blue,
    filecolor=blue,      
    urlcolor=blue,
    citecolor=blue,
    }

\usepackage[mathlines]{lineno}

\usepackage[margin=1in]{geometry}

\numberwithin{equation}{section}

\def\d{\,\mathrm{d}}
\def\N{\mathbb{N}}
\def\R{\mathbb{R}}
\def\C{\hbox{\rlap{\kern.24em\raise.1ex\hbox
      {\vrule height1.3ex width.9pt}}C}}
\def\P{\hbox{\rlap{I}\kern.16em P}}
\def\Q{\hbox{\rlap{\kern.24em\raise.1ex\hbox
      {\vrule height1.3ex width.9pt}}Q}}
\def\M{\hbox{\rlap{I}\kern.16em\rlap{I}M}}
\def\Z{\hbox{\rlap{Z}\kern.20em Z}}
\def\({\begin{eqnarray}}
\def\){\end{eqnarray}}
\def\[{\begin{eqnarray*}}
\def\]{\end{eqnarray*}}
\def\part#1#2{\frac{\partial #1}{\partial #2}}

\def\grad{\nabla}

\def\bar{\overline}

\def\tot#1#2{\frac{\d #1}{\d #2}} 

\def\laplace{\Delta}
\def\d{\,\mathrm{d}}
\def\N{\mathbb{N}}
\def\R{\mathbb{R}}

\def\epsilon{\varepsilon}

\def\d{\mathrm{d}}

\usepackage{authblk}

\def\rev#1{#1}
\def\revv#1{#1}

\title{Impact of memory on clustering in spontaneous particle aggregation}
\date{}

\author{
     Radek Erban\\ Mathematical Institute, University of Oxford,\\
     Radcliffe Observatory Quarter, Woodstock Road,\\
     Oxford, OX2 6GG, United Kingdom\\
     {\it radek.erban@maths.ox.ac.uk}\\
    \and
     Jan Haskovec\\ Mathematical and Computer Sciences
            and Engineering Division,\\
         King Abdullah University of Science and Technology,\\
         Thuwal 23955-6900, Kingdom of Saudi Arabia\\
         {\it jan.haskovec@kaust.edu.sa}
     }

\begin{document}

\maketitle

\begin{abstract} \noindent
The effect of short-term and long-term memory on spontaneous aggregation of
organisms is investigated using a stochastic agent-based model.
Each individual modulates the amplitude of its random motion according to the perceived local density of neighbors.
Memory is introduced via a chain of $K$~internal  variables
that allow agents to retain information about previously encountered densities.
The parameter $K$ controls the effective length of memory.
A formal mean-field limit yields a macroscopic Fokker--Planck equation, which provides a continuum description of the system in the large-population limit.
Steady states of this equation are characterized to interpret the emergence and morphology of clusters.
Systematic stochastic simulations in one- and two-dimensional spatial domains reveal that short- or moderate-term memory promotes coarsening,
resulting in a smaller number of larger clusters, whereas long-term memory inhibits aggregation and increases the proportion of isolated individuals. 
Statistical analysis demonstrates that extended memory reduces the agents' responsiveness to environmental stimuli,
explaining the transition from aggregation to dispersion as $K$ increases. 
These findings identify memory as a key factor controlling the collective organization of self-driven agents
and provide a bridge between individual-level dynamics and emergent spatial patterns.

\textbf{Relevance to Life Sciences.}
The presented framework provides a quantitative approach for studying how biological memory influences collective organization. Many biological systems,
including bacterial populations, cellular assemblies, and insect groups, display aggregation behavior that depends on how past environmental information is integrated into current decision-making.
The model bridges individual cognition (through memory) and population-level spatial structure (through clustering).
It demonstrates that limited memory enhances aggregation by maintaining sensitivity to current local conditions, while excessively long memory diminishes responsiveness and disrupts coordination. This mechanism offers a potential explanation for experimentally observed adaptation timescales in biochemical and behavioral systems, such as the short-term biochemical memory of bacterial chemotaxis or transient social memory in gregarious insects. The results suggest that the duration of memory may have evolved to balance responsiveness and stability in collective behaviors, leading to efficient spatial organization under fluctuating environmental conditions.

\textbf{Mathematical Content.}
The model is formulated as a system of coupled stochastic differential equations for $N$ interacting agents,
each endowed with $K$ internal variables.
A system of delayed stochastic integro-differential equations is derived with its memory kernel quantifying how past
information influences particle motion and determines the effective memory timescale. 
The formal mean-field limit, derived under the molecular chaos assumption, yields a nonlinear and nonlocal Fokker--Planck equation
for the joint density of spatial and internal states. Analytical characterization of its steady states identifies conditions
for spatially inhomogeneous equilibria corresponding to particle clusters.
Numerical investigations employ direct stochastic simulations combined with density-based spatial clustering analysis
to quantify aggregation statistics and elucidate how memory influences agents' responsiveness to environmental stimuli.
The study thus links microscopic stochastic processes with macroscopic pattern formation through a unified mathematical framework.
\end{abstract}

\vskip 5mm


\vskip 5mm


\section{Introduction}
\label{sec1}

\noindent
Models with memory and delays are ubiquitous in explaining collective behaviour of organisms,
ranging in scale from bacterial and amoeboid chemotaxis~\cite{Berg:1972:CEC,Erban:2007:TEA} to behaviour 
of social insects and flocks of birds~\cite{Sumpter:2006:PCA,Sumpter:2010:CAB}. In the case of bacteria, 
the memory is incorporated in purely biochemical terms, in the form of the signal transduction network 
capable of excitation and adaptation dynamics~\cite{Spiro:1997:MEA,Alon:1999:RBC}. This enables a cell 
to `memorize' past environmental signals and compare to their current state to inform decision 
making~\cite{Berg:1990:BM}. In individual-based models, such memory effects can be accounted for
by ordinary differential equation (ODE) or stochastic differential equation (SDE) models describing 
internal variables of each individual~\cite{Erban:2004:ICB,Erban:2007:TEA}.
A typical example of such internal variables are concentrations of key intracellular biochemical species.
At the collective (macroscopic) level, populations of bacteria or other cells can be described by partial 
differential equations (PDEs) describing the density of cells and key environmental signals~\cite{Erban:2005:STS,Erban:2020:SMR}, 
where memory can also be interpreted as delays in the processing of some signals. Incorporating such delays 
into mathematical models of collective phenomena can provide more accurate explanations of behavioural 
properties of groups of animals~\cite{Sun:2014:TDC,Erban:2016:CSM}, robots~\cite{Taylor:2015:MMT},
or spreading processes in social and biological systems~\cite{Rosvall:2014}.
Moreover, non-Markovian characteristics resulting from ordering of interactions in temporal complex networks
were identified as an important mechanism that alters causality and affects dynamical processes
in social, technological and biological systems~\cite{Scholtes2014, Williams2022}.
Considering biological agents (cells or animals), a trade-off between having no memory at all and remembering too 
much (including very distant and irrelevant past states) has been achieved by evolution.
For example, bacterial `biochemical memory' of past environmental signals evolved to last for the duration 
of a few seconds~\cite{Berg:1990:BM}.

In this paper we focus on a model of spontaneous aggregation of animals or cells resulting
from random diffusive motion of individuals~\cite{BHW:2012:PhysD, HO:2015:DCDSB}.
The individuals respond to the local population density observed in their neighbourhood
by increasing or decreasing the amplitude of their random motion. This kind of behaviour 
has been observed in insects, for instance, the pre-social German cockroach (\textit{B. germanica}),
which are known to be attracted to dark, warm and humid places. However, 
it has been shown that cockroach larvae also aggregate spontaneously~\cite{Jeanson:2005:AnimBehav}, 
{i.e.}, in the absence of any environmental template or heterogeneity.
Moreover, this type of dynamics can also be used to describe formation of P Granules
in \textit{C. elegans} embryos resulting in spontaneous protein aggregation~\cite{Brangwynne:2009:Science}.
\rev{We note that spontaneous aggregates form at random locations in space.  If such spontaneous aggregates 
are to form on a regular lattice, additional interactions between individuals must be incorporated into the 
model~\cite{Lopez:2025:SCR}.}

To incorporate memory into the first-order model of spontaneous aggregation~\cite{BHW:2012:PhysD},
we introduce a set of $K\geq 1$ internal variables of each agent. \rev{Such internal variables can have
different physical meanings in different biological models. In applications to bacterial chemotaxis, individuals
compare current and past concentrations of extracellular signals, and internal variables can be identified
with the concentration of molecules transducing the signal in detailed models of the chemotactic network~\cite{Celani:2010}.
The internal variables can also be viewed as effectively storing the memory of the state of the external
environment in the vicinity of the individual. In our case, we consider a model of the spontaneous aggregation 
(i.e., with no external chemical signal) and $K \ge 1$ internal variables store a time-dependent memory of how 
crowded the neighborhood was in the past. In particular, each internal variable} describes 
a `layer' of memory, which is subject to two effects: (i)~production (or excitation) from the internal 
variable of the next higher order, and (ii)~spontaneous decay with a constant rate. The internal variable 
of the highest ({i.e.}, $K$-th) order is then subject to random excitation, with amplitude modulated by 
a nonlinear function of the number density of agents observed in the physical neighbourhood. We note 
that in the first-order model introduced in~\cite{BHW:2012:PhysD}, the positions of the agents 
in the physical space were directly subjected to random excitation (Brownian motion). In particular,
our investigation extends their model by introducing a chain of $K$ internal variables that allow the 
agents to `remember' the densities they encountered in the past.

Our results show that the introduction of memory with a few layers, $K$, leads to the formation of a 
smaller number of larger clusters of agents. This trend is observed until $K=3$ or $K=4$, depending 
on the spatial  dimension of the studied system, as shown in Sections~\ref{sec5} and~\ref{sec6}. 
When the number of memory layers $K$ is increased further, {i.e.}, as the memory becomes `longer', 
its effect starts to be disruptive. This is manifested by the increasing proportion of `outliers', which 
are the agents that are not part of any cluster. We therefore conclude that short-term or medium-term 
memory enhances spontaneous aggregation, while long-term memory disrupts it.
\rev{
This observation is ecologically relevant because many organisms - from bacteria to insects to higher animals -
do not respond to their environment instantaneously but instead use some form of memory or internal state to modulate their behavior.
Species that benefit from aggregation (e.g., for warmth, protection, or reproduction)
may have developed memory systems that retain only short-to-medium-term environmental information.
Conversely, overly long memory may reduce responsiveness to current cues and hinder effective social or spatial organization.
Therefore, organisms must balance between being too reactive (leading to unstable or erratic behavior)
and being too inertial (failing to adapt to rapidly changing environments).
This is consistent with ecological and evolutionary studies that view memory
as an adaptive trait shaped by natural selection, see, e.g.,~\cite{Couzin:2002, Chen:2024}.
}

\rev{
There is a vast body of literature focusing on mathematical models of aggregation dynamics and their analysis,
with the number of papers counting in hundreds. These models
include not only aggregation of agents or particles in the physical space
(e.g.,~\cite{Jeanson:2005:AnimBehav, Brangwynne:2009:Science, Lopez:2025:SCR, Keller:1970:JTB, Haskovec:2009:JSP}),
but also in abstract phase space, i.e., models of
opinion formation~\cite{Hegselmann:2002:JASSS, Jabin:2014:JDE},
flocking and swarming~\cite{Vicsek:1995:PRL, Cucker:2007:IEEE, Olfati:2006:IEEE}
or pattern formation resulting from attraction-repulsion interactions~\cite{Dorsogna:2006:PRL, Carrillo:2019:AP}.
For a  relatively recent survey of models of collective motion we refer to~\cite{Vicsek:2012:PR}.
Many variants of these models have been studied, including models with nonlocal interactions~\cite{Mogilner:1999:JMB, Bernoff:2013:SIAMRev, Carrillo:2015:KRM},
topological interactions~\cite{Haskovec:2013:PhysD}, hydrodynamic models~\cite{Choi:2019:SIMA}, models with hierarchy~\cite{Shen:2007:SIAM, Voggt:2019:PhysRevRes} and leaders~\cite{Cicolani:2024:arXiv},
models with noise~\cite{Cattiaux:2018:AAP} and delay~\cite{Choi:2021:MMAS, Haskovec:2021:SIADS} and combination of those~\cite{Erban:2016:CSM}.
In~\cite{Ahn:2023:MBE} the authors considered a variant of the Cucker-Smale flocking model
with a Caputo fractional derivative in the time variable. The presence of the fractional derivative
can be interpreted as a form of memory. The authors derived sufficient conditions for asymptotic flocking,
which occurs at an algebraic rate, in contrast to the exponential flocking asymptotics in the original Cucker-Smale model~\cite{Cucker:2007:IEEE}.
Memory induced aggregation in collective foraging was studied in~\cite{Nauta:2020:HFP}.
Using an agent-based model, the authors argued that aggregation around salient patches
can occur through formation of collective memory realized through local interactions and global
displacement using L\'{e}vy walks.
The impact of emergent memory on opinion dynamics was analyzed in~\cite{Boschi:2021:PhysA},
with the aim to explain how the exposure of the society to certain events deeply
changes people's perception of the present and future.
An analytical way was proposed to measure how much information a society
can remember when an extensive number of news items was presented.
}

\rev{
In all the aforementioned models, aggregation arises from explicit deterministic forces acting on the particles,
typically expressed as a superposition of pairwise interactions.
Eventual random effects only disrupt the aggregation.
In contrast, in the spontaneous aggregation model studied in this paper,
each particle is merely subject to Brownian motion
with amplitude modulated by a nonlinear function of the density.
In particular, on the microscopic level of description, there is no deterministic force acting on the particles
and the interaction between them is indirect. 
In particular, the model cannot be formulated in terms of pairwise particle interactions.
Therefore, aggregation in our model arises solely from a stochastic process
and we cannot expect formation of regular patterns, like flocks, swarms or opinion clusters
as in the aforementioned models.
}

The paper is organized as follows. In Section~\ref{sec2} we describe the (first-order) spontaneous 
aggregation model and summarize its main properties studied in the literature~\cite{BHW:2012:PhysD}. 
In Section~\ref{sec3} we introduce memory into this model and infer the main mathematical properties
of the model with memory. In Section~\ref{sec:FP} we derive the formal macroscopic limit of the system 
as the number of agents tends to infinity, obtaining the corresponding Fokker-Planck equation.
In Section~\ref{sec:patterns} we then characterize
its steady states to gain an insight into the patterns (clusters) formed by the system. 
In Section~\ref{sec5} we report the results of the stochastic simulations of the individual-based model 
in the spatially one-dimensional setting, while the two-dimensional results are presented in Section~\ref{sec6}. 
Here we also provide a statistical evidence that the long-term memory inhibits the particles' responsiveness 
to environmental stimuli. We conclude with the discussion in Section~\ref{sec7}.

\section{Spontaneous aggregation model without memory}

\label{sec2}

\noindent
The individual-based stochastic model introduced in reference~\cite{BHW:2012:PhysD}
under the name `direct aggregation model' consists of a group of $N \ge 2$ biological agents (cells or animals),
characterized by their positions ${\mathbf x}_i(t)\in\R^d$, with spatial dimension $d \in \{1,2,3\}$ and $i\in [N]$,
where we have denoted the set of indices by $[N] := \{1,2,\ldots,N \}$. Every individual senses the average 
density of its close neighbours, given by
\begin{equation}   
\label{rho_i}
\vartheta_i(t) = \frac{1}{N} \sum_{j\in [N]} W({\mathbf x}_i(t)-\mathbf{x}_j(t)) \,,\qquad \mbox{for} \quad i \in [N]\,,
\end{equation}
where $W({\mathbf x}) = w(|{\mathbf x}|)$ with the weight function $w: \R^+ \to \R^+$ assumed to be bounded, nonnegative, nonincreasing and
integrable on $\R^d$. Without loss of generality we impose the normalization
\begin{equation}  
\label{eq:Wnorm}
\int_{\R^d} W({\mathbf x}) \, \d {\mathbf x} = 1 \,. 
\end{equation}
A generic example of $w$ is the (properly normalized) characteristic function
of the interval $[0,R]$, corresponding to the sampling radius $R > 0$.
The average density $\vartheta_i$ is then simply the fraction
of individuals located within the distance $R$ from the $i$-th individual.
The individual positions are subject to a random walk with modulated amplitude,
described by the system of coupled SDEs
\begin{equation}
\label{model1}
\d {\mathbf x}_i(t) \,=\, G(\vartheta_i) \, \d {\mathbf B}_i^t \,,\qquad \mbox{for} \quad i \in [N]\,,
\end{equation}
where ${\mathbf B}_i^t$ are independent $d$-dimensional Brownian motions.
The \emph{response function} $G:\R^+ \to\R^+$ is assumed to be globally bounded, nonnegative and decreasing.
The monotonicity of $G$ is implied by the modeling assumption that the individuals respond to 
higher perceived population densities in their vicinity by reducing the amplitude of 
their random walk.

In reference~\cite{BHW:2012:PhysD}, the mean-field limit (as $N\to\infty$) of the individual-based model~\eqref{rho_i}--\eqref{model1}
is studied, describing the system in terms of the number particle density $\varrho \equiv \varrho(t,{\mathbf x})$
with ${\mathbf x}\in\R^d$. The time evolution of the density $\varrho$ is subject to
\begin{equation}   
\label{eq:rho}
\part{\varrho}{t} = \frac12 \, \laplace\left( G(W\ast\varrho)^2 \varrho \right),
\end{equation}
with the convolution given by~$W\ast\varrho({\mathbf x}) := \int_{\R^d} W({\mathbf x}-{\mathbf z}) \varrho({\mathbf z}) \, \d {\mathbf z}$.
In fact, it is more instructive to expand the derivative in equation~\eqref{eq:rho}
and formulate it as the convection-diffusion equation
\begin{equation}   
\label{eq:convdiff}
\part{\varrho}{t} =  \frac12 \, \grad\cdot\Big( \varrho \, \grad\big[ G(W\ast\varrho)^2\big] + G(W\ast\varrho)^2 \, \grad \varrho \Big).
\end{equation}
Here we identify the convection term $\grad\cdot\left( \varrho \, \grad\big[ G(W\ast\varrho)^2\big]\right)$ 
as the `driving force' responsible for the eventual formation of aggregates
(recall the monotonicity assumption on $G$). Then the following condition for aggregation, emanating from a perturbation
of a given constant steady state \mbox{$\varrho \equiv \varrho_0 > 0$}, can be derived~\cite{BHW:2012:PhysD},
\begin{equation}   
\label{cond:agg}
\mbox{Re } \hat W(\boldsymbol{\xi}) \,>\, - \frac{G(\varrho_0)}{2\,G'(\varrho_0)\,\varrho_0} \quad\qquad\mbox{for some } \boldsymbol{\xi}\in\R^d,
\end{equation}
where $\hat W = \hat W(\boldsymbol{\xi})$ is the Fourier transform of $W$
and $G'(\varrho_0)$ denotes the derivative of $G$ at $\varrho_0$.
Note that due to the monotonicity assumption, we have $G'(\varrho_0) < 0$.
For a fixed kernel $W=W({\mathbf x})$ one can interpret the inequality~\eqref{cond:agg} as
a condition on the response function $G$ to be steep enough
in the neighbourhood of $\varrho_0$.
Then the effect of the convection term in equation~\eqref{eq:convdiff} is stronger
than the smoothing effect of the diffusion, and the resulting instability
leads to formation of aggregates. These aggregates persist as a steady state solution
in the large time limit.

\section{Spontaneous aggregation model with memory}

\label{sec3}

\noindent
We now introduce memory into the spontaneous aggregation model~\eqref{rho_i}--\eqref{model1}
using a set of $K \in \N$ internal variables describing $K$ memory `layers'. Each individual is 
then characterised by its position ${\mathbf x}_i={\mathbf x}_i(t)\in\R^d$
and internal variables ${\mathbf y}_{i}^{k}={\mathbf y}_{i}^{k}(t)\in\R^d$ for $i \in [N]$ and $k\in [K]$.
Equation~(\ref{model1}) is generalized for $i \in [N]$ to the following system of SDEs,
\begin{equation}   
\label{SDE1}
\begin{aligned}
\tot{{\mathbf x}_i(t)}{t} &\,=\, {\mathbf y}_{i}^1(t)\,, \\
\varepsilon_k \tot{{\mathbf y}_{i}^k(t)}{t} &\,=\, - \alpha_{k} \,{\mathbf y}_{i}^k(t) \,+\, {\mathbf y}_{i}^{k+1}(t)\,, \qquad\quad \mbox{for} \;\; k \in [K-1]\,, \\
\varepsilon_K \,\d {\mathbf y}_{i}^K(t) &\,=\,  - \alpha_K \,{\mathbf y}_{i}^K(t) \, \d t \,+\, G(\vartheta_i(t)) \,\d {\mathbf B}_i^t\,,
\end{aligned}
\end{equation}
where ${\mathbf B}_i^t$ are again independent $d$-dimensional Brownian motions.
The positive constants $\alpha_k > 0$, for $k\in[K]$, are the relaxation coefficients,
and parameters $\varepsilon_k>0$, for $k\in[K]$, define the time scales.
The perceived densities $\vartheta_i = \vartheta_i(t)$ are given by \eqref{rho_i}.

\revv{The internal variables ${\mathbf y}_{i}^{k}={\mathbf y}_{i}^{k}(t)\in\R^d$
represent successive layers of biological memory through which information about past local crowding
is processed before affecting movement. Rather than storing the past density explicitly,
the agent integrates past sensory input through a cascade of internal states with finite relaxation times.
Such cascades naturally arise in biological systems, for instance in biochemical signalling networks
where an external stimulus is transmitted through multiple molecular intermediates,
or in neural or behavioral processing where information is gradually filtered and attenuated.
}

\revv{
In this interpretation, the parameter $K$ does not represent distinct memories of specific past times,
but rather the depth of the memory-processing cascade.
Increasing $K$ increases the number of intermediate steps through which information must pass before influencing motion,
thereby extending the effective memory timescale and smoothing the response to rapid environmental fluctuations.
}

\rev{The model with memory~(\ref{SDE1}) has been introduced in a way that it can be reduced to the original 
spontaneous aggregation model~\eqref{rho_i}--\eqref{model1} if we set $\alpha_k:=1$ for all 
$k\in[K]$ and pass to the limit $\varepsilon_k \to 0$ for all $k\in[K]$. However, it can be also formulated
in a more general way, which provides connections with the memory models in the literature~\cite{Fagan:2023:RDM, Nauta:2020:HFP, Kim:2024:ISM}.}
Denoting the $j$-th component of ${\mathbf y}_{i}^{k} \in\R^d$ by $y_{i,j}^k$, for 
$i \in [N]$, $j\in [d]$ and $k\in [K]$, we introduce the following notation
\begin{equation}
{\mathbf y}_{i,j} = \big[y_{i,j}^1, y_{i,j}^2, \dots, y_{i,j}^K\big] \in \R^{K}
\qquad\quad \mbox{for} \quad i\in[N] \;\; \mbox{and} \;\;j\in [d]\,.
\label{internalvariablenotation}
\end{equation} 
To further analyze the behaviour of model~(\ref{SDE1}), we fix 
\begin{equation}
\varepsilon_k:=1 \qquad\quad \mbox{for all} \quad k\in[K].
\label{epsare1} 
\end{equation}
Then the SDE system for ${\mathbf y}_{i,j} = {\mathbf y}_{i,j}(t) \in \R^K$,
with $i\in[N]$ and $j\in [d]$, can be written in matrix form as
\begin{equation}
\label{SDE}
\d{{\mathbf y}_{i,j}(t)} \,=\, A\rule{0pt}{9pt}^{(K)}\, {\mathbf y}_{i,j}(t) \, \d t \,+\, {\boldsymbol{\beta}}^{(K)}_i(t) \, \d B_{i,j}^t
\end{equation}
where $B_{i,j}^t$ is the $j$-th component of ${\mathbf B}_{i}^t$ ({i.e.}, the standard one-dimensional Brownian motion), and
the constant matrix $A\rule{0pt}{9pt}^{(K)} \in \R^{K\times K}$ and the vectors ${\boldsymbol{\beta}}^{(K)}_i = {\boldsymbol{\beta}}^{(K)}_i(t) \in \R^{K}$
are given by
\begin{equation}
\label{AKBK}
A\rule{0pt}{9pt}^{(K)} := \begin{pmatrix}
      - \alpha_1 & 1 & 0 & 0 & \ldots & 0 \\
      0 & -\alpha_2 & 1 & 0 & \ldots & 0 \\
      \vdots & & \ddots & & 1 & 0 \\
      0 & &  & & -\alpha_{K-1}& 1 \\
      0 & & \ldots & & 0 & -\alpha_K
   \end{pmatrix}
\qquad\quad\mbox{and}\qquad\quad
{\boldsymbol{\beta}}^{(K)}_i(t) := \begin{bmatrix} 0 \\ 0 \\ \vdots \\ 0 \\ G(\vartheta_i(t)) \end{bmatrix}.
\end{equation}
\rev{The reformulation of model~(\ref{SDE1}) into the SDE form~(\ref{SDE}) aligns it more closely with the formulation of the reinforced diffusion model for memory-mediated animal movement~\cite{Fagan:2023:RDM}. The key difference in our study is that we focus on spontaneous aggregation, whereas some models in the literature examine aggregation in response to environmental cues to locate specific targets~\cite{Nauta:2020:HFP, Kim:2024:ISM}. Additionally, we explicitly describe the memory variables, whereas some models introduce memory implicitly through delays~\cite{Fagan:2023:RDM}. Model~(\ref{SDE}) can also be rewritten in the delay form, because the solution of the SDE system~\eqref{SDE} is given} by the variation-of-constants formula~\cite[Section 3.3]{Mao:2007}
\begin{equation}
{\mathbf y}_{i,j}(t) 
= \exp\!\Big[t A\rule{0pt}{9pt}^{(K)}\Big] \, {\mathbf y}_{i,j}(0) \,+\,  \int_0^t \exp\!\Big[(t-s) A\rule{0pt}{9pt}^{(K)}\Big] \, {\boldsymbol{\beta}}^{(K)}_i(s) \, \d B_{i,j}^s \,,
\label{varformula}
\end{equation}
where ${\mathbf y}_{i,j}(0)$ is a prescribed initial condition at time $t=0$. \rev{The SDE form of the memory model~(\ref{SDE}) can, in general, be used to describe the `chemical memory' of unicellular organisms, where internal variables ${\mathbf y}_{i}^{k}$ can be identified with concentrations of key chemical species of signal transduction networks~\cite{Erban:2005:STS,Erban:2007:TEA}. The behaviour of such memory models will depend on the properties of matrix $A^{(K)} \in \R^{K\times K}$, on its eigenvalues and eigenvectors. In this paper, we use the specific form of the general model~(\ref{SDE}) given by~(\ref{AKBK}).} In particular, choosing ${\mathbf y}_{i,j}(0) = 0$ in equation~(\ref{varformula}), we have for $y_{i,j}^1 \equiv y_{i,j}^1(t)$
$$
y^1_{i,j}(t) = \int_0^t \kappa(t-s) \,G(\vartheta_i(s)) \,\d B_{i,j}^s,
\qquad
\mbox{with} 
\qquad
\kappa(t) := \left\{ \exp\!\Big[t A\rule{0pt}{9pt}^{(K)}\Big] \right\}_{1,K}, 
$$
where $\kappa(t)$ is the $(1,K)$-th element of the matrix $\exp\!\Big[t A^{(K)}\Big] \in \R^{K\times K}$. To get further insight
into the model behaviour we assume
\begin{equation}
\alpha_k:=\alpha \qquad\quad \mbox{for all} \quad k\in[K]\,,
\label{alphakarealpha} 
\end{equation}
where $\alpha>0$ is a single parameter of the model. Then the kernel $\kappa=\kappa(t)$ can be evaluated explicitly,
because the matrix $A^{(K)}$ in \eqref{AKBK} takes the form
$$
A\rule{0pt}{9pt}^{(K)} \,=\, -\alpha \,I\, +\, U \,,
$$
with $U$ the upper-diagonal matrix with ones on the first upper diagonal. \rev{In particular, all eigenvalues of matrix $A^{(K)}$ are equal to $-\alpha$.} Since $U$ commutes with the identity matrix, we have
\[
\exp\!\Big[t A\rule{0pt}{9pt}^{(K)}\Big]  \, = \, \exp\left( -\alpha t \right) \,\exp[t U] \,.
\]
Moreover, $U$ being upper diagonal, its $K$-th power $U^K$ is the zero matrix, so that
\[
   \exp\left(t U \right) = \sum_{k=0}^{K-1} \frac{t^k U^k}{k!}.
\]
Consequently, the $(1,K)$-th element of $\exp\left[t U \right]$ equals to $t^{K-1}/(K-1)!$
and we have
\(  \label{eq:kappa}
\kappa(t) = \left\{ \exp\!\Big[t A\rule{0pt}{9pt}^{(K)}\Big] \right\}_{1,K} = \frac{e^{-\alpha t} \, t^{K-1}}{(K-1)!}\,.
\)
Using equation~\eqref{SDE1}, we obtain
\begin{equation}   
\label{eq:SDE2}
\tot{{\mathbf x}_i(t)}{t} = \frac{1}{(K-1)!} \int_0^t e^{-(t-s)\alpha} \, (t-s)^{K-1} \, G(\vartheta_i(s)) \,\d {\mathbf B}_{i}^s\,,
\end{equation}
with $\vartheta_i = \vartheta_i(s)\in\R$ given by equation~\eqref{rho_i}.
Let us note that this formulation allows for generalization of the model
to noninteger values of $K>0$ by writing the factorial $(K-1)!$
in terms of the Gamma-function as $\Gamma(K)$.

\revv{
To connect the mathematical formulation with a biologically meaningful notion of memory duration,
we characterize the effective memory length as the 
mean past time interval that contributes to the current motion of an individual.
This quantity corresponds to the mean of the normalized version of the kernel $\kappa=\kappa(t)$,
}
\( \label{eq:length}
\left( \int_0^{+\infty} \kappa(t) \, \d t \right)^{\!\!-1} \left( \int_0^{+\infty} t \, \kappa(t) \, \d t \right) = \frac{K}{\alpha}.
\)
 \revv{
 The quantity $K/\alpha$ therefore provides a direct measure of the biological memory timescale.
 It represents the average time period between when local crowding is perceived
 and when it influences the agent's movement. Increasing the number of memory layers 
$K$ proportionally increases this period, meaning that individuals respond to a longer history
of past environmental conditions rather than primarily to the recent past.
From a biological perspective, larger values of $K$ correspond to organisms with slower adaptation
or longer-lasting internal states, while smaller $K$ corresponds to short-term memory dominated by recent sensory input.
 }
 
\begin{figure}
\begin{center}
\includegraphics[width=.45\textwidth]{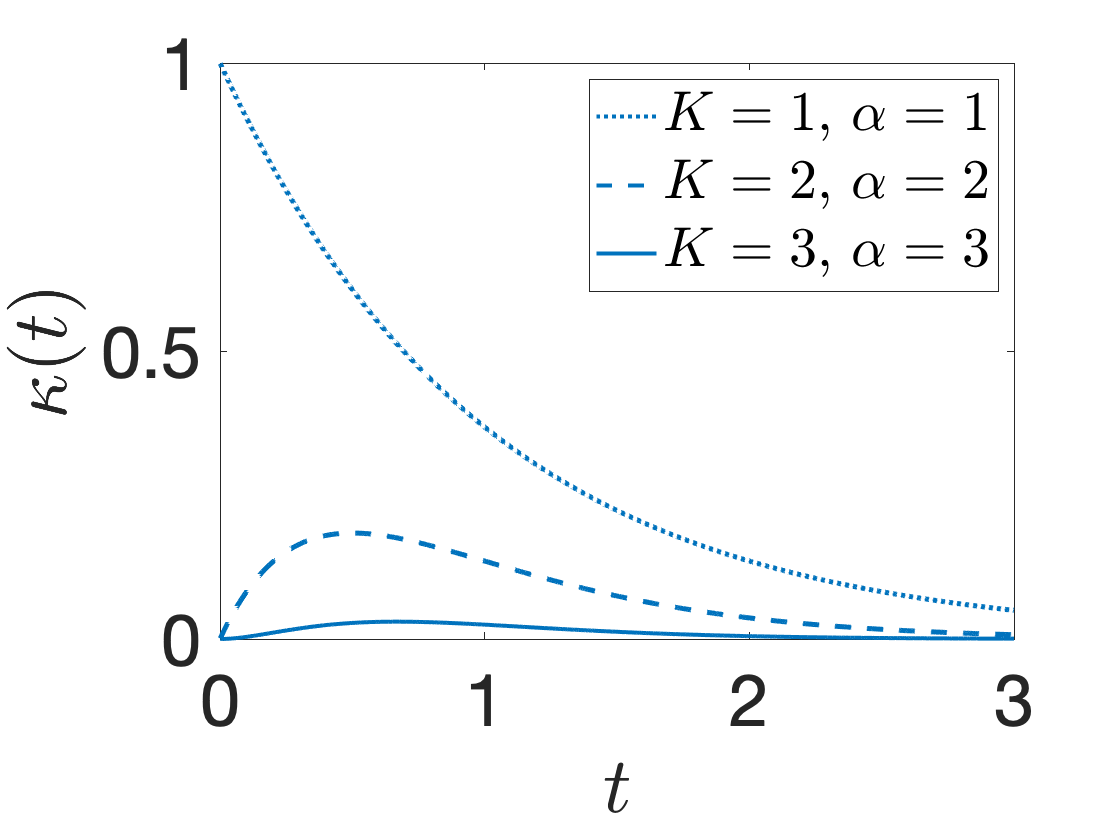}
\includegraphics[width=.45\textwidth]{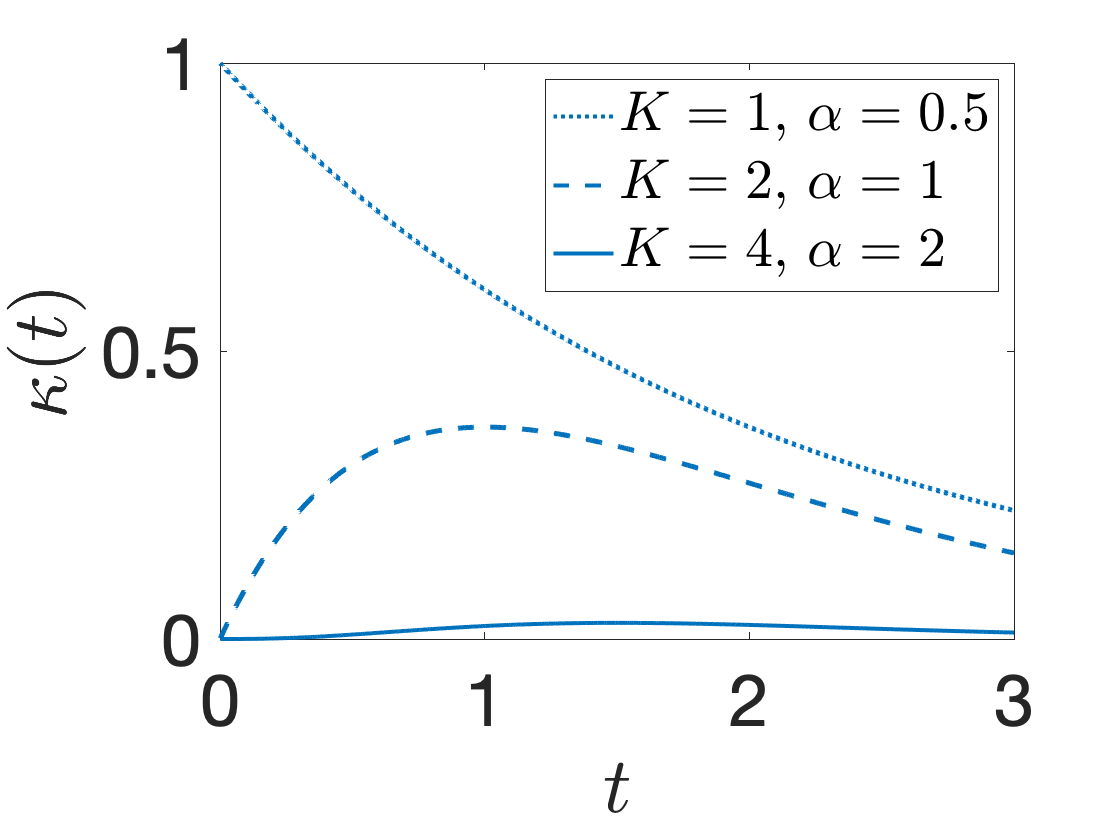}
\end{center}
\caption{{\it Plots of the kernel $\kappa=\kappa(t)$ given by \eqref{eq:kappa} with $K/\alpha = 1$ (left panel)
and with $K/\alpha=2$ (right panel).
}}
\label{fig:kappa}
\end{figure}

 In Sections~\ref{sec5} and~\ref{sec6}, 
we will numerically investigate the impact of the choice of $K$ on the clustering properties of the individual-based model 
with memory given by equations~\eqref{rho_i} and~\eqref{eq:SDE2}.
\rev{
Let us note that the kernel \eqref{eq:kappa} takes very different shapes
even if the length \eqref{eq:length} is kept constant.
For instance, in Figure~\ref{fig:kappa} we plot $\kappa=\kappa(t)$ for the length 
$$\frac{K}{\alpha} = \left\{ \frac11, \frac22, \frac33 \right\} = 1 \mbox{ (left panel),}
\qquad \mbox{and} \qquad
\frac{K}{\alpha} = \left\{ \frac1{0.5}, \frac21, \frac42 \right\} = 2 \mbox{ (right panel).}
$$
In particular, for $K=1$ the kernel reduces to $\kappa(t) = e^{-\alpha t}$, i.e., is globally strictly decreasing on $[0,+\infty)$,
while for $K\geq 2$ it has a unique strict maximum in $(0,+\infty)$.
Therefore, it is not sufficient to only explore the case $K=1$ and tune the value of $\alpha>0$
to modulate the length of the memory~\eqref{eq:length}.
In this paper we chose to fix the value $\alpha=1$ and explore the impact of different choices
of the number of layers $K$ on the clustering properties of the system.
 }%

Finally, let us note that this individual-based model 
can be also viewed  for $K=1$ as a simplified version of what is called `the second-order model' in reference~\cite{BHW:2012:PhysD}. 
Indeed, if $K=1$, then the SDE system~\eqref{SDE1} reduces to
\begin{equation}
\label{SDExv}
\begin{aligned}
\tot{{\mathbf x}_i(t)}{t} &\,=\, {\mathbf y}_{i}^1(t)\,, \\
\d {\mathbf y}_{i}^1(t) &\,=\, - \alpha \, {\mathbf y}_{i}^1(t) \,\d t \,+\, G(\vartheta_i(t)) \,\d {\mathbf B}_i^t\,.
\end{aligned}
\end{equation}
In particular, we can interpret the internal variable ${\mathbf y}_{i}^1 = {\mathbf y}_{i}^1(t)$
as the velocity of the $i$-th agent in the physical space.
The velocity is then subject to linear damping with intensity $\alpha>0$
and stochastic forcing with amplitude modulated by the term $G(\vartheta_i(t))$
with $\vartheta_i = \vartheta_i(t)$ given by~\eqref{rho_i}.

\section{Formal derivation of the Fokker-Planck equation}
\label{sec:FP}

\noindent
The formal large population limit as $N\to\infty$ of the system~\eqref{SDE1} with the average density given by~\eqref{rho_i}
can be derived by a generalization of the procedure used in~\cite[Section 3]{BHW:2012:PhysD}.
Since some nontrivial steps were omitted in~\cite{BHW:2012:PhysD}, we include the detailed formal derivation here.
The main idea is to `linearize'~\eqref{SDE1} by introducing the average densities
$\vartheta_i = \vartheta_i(t)$ for $i\in[N]$ as a set of $N$~independent variables. According 
to~\eqref{rho_i} and the first equation of~\eqref{SDE1}, these new variables
are subject to the system of differential equations
$$
\tot{\vartheta_i}{t} \,=\, \frac{1}{N} \sum_{j\in[N]} \grad W({\mathbf x}_i - {\mathbf x}_j) \cdot \left({\mathbf y}_i^1 - {\mathbf y}_j^1 \right), 
\qquad\quad \mbox{for} \;\; i\in [N]\,.
$$
One then applies the It\^{o} formula~\cite{Mao:2007} to write down the Liouville equation~\cite{Cercignani:1988}
for the $N$-particle density
$f^N = f^N(t, {\mathbf x}_1, \dots, {\mathbf x}_N, {\mathbf y}^1_1, \dots, {\mathbf y}^1_N, \dots,
{\mathbf y}^K_1, \dots, {\mathbf y}^K_N, \vartheta_1, \dots, \vartheta_N)$,
\[
   \part{f^N}{t} + \sum_{i\in[N]} {\mathbf y}^1_i \cdot\grad_{{\mathbf x}_i} f^N
   + \sum_{i\in[N]} \sum_{k\in[K-1]} \frac{1}{\varepsilon_k} \grad_{{\mathbf y}^k_i} \cdot \left( ({\mathbf y}^{k+1}_i - \alpha_k {\mathbf y}^{k}_i) f^N \right) \\
   - \frac{1}{\varepsilon_K} \sum_{i\in[N]} \grad_{{\mathbf y}^K_i} \cdot \left( \alpha_K {\mathbf y}^{K}_i f^N \right) \\
   + \frac1N \sum_{i\in[N]} \sum_{j\in[N]} \grad W({\mathbf x}_i - {\mathbf x}_j) \cdot \left({\mathbf y}_i^1 - {\mathbf y}_j^1 \right) \part{f^N}{\vartheta_i} \\
   = \frac{1}{2\varepsilon_K} \sum_{i\in[N]} \laplace_{{\mathbf y}_i^K} \left( G(\vartheta_i)^2 f^N \right).
\]
Integrating with respect to the variables $({\mathbf x}_i, {\mathbf y}^1_i, \dots, {\mathbf y}_i^K, \vartheta_i)$ for $i=2,\dots,N$
and adopting the usual molecular chaos assumption on vanishing particle correlations,
we obtain
the Vlasov equation for the one-particle marginal density $f^1 = f^1(t, {\mathbf x}, {\mathbf y}^1, {\mathbf y}^2, \dots, {\mathbf y}^K, \vartheta)$.
For brevity, we introduce the notation ${\mathbf y} := ({\mathbf y}^1, {\mathbf y}^2, \dots, {\mathbf y}^K)$. Then, the Vlasov equation
for $f^1 = f^1(t, {\mathbf x}, {\mathbf y}, \vartheta)$ reads
\[
   \part{f^1}{t} + {\mathbf y}^1 \cdot\grad_{\mathbf x} f^1
   +  \sum_{k=1}^{K-1} \frac{1}{\varepsilon_k} \grad_{{\mathbf y}^k} \cdot \left( ({\mathbf y}^{k+1} - \alpha_k {\mathbf y}^{k}) f^1 \right) \\
   - \frac{1}{\varepsilon_K} \grad_{{\mathbf y}^K} \cdot \left( \alpha_K {\mathbf y}^{K} f^1 \right) 
   + \part{}{\vartheta} 
   \left( F[f^1] f^1 \right) \\
   = \frac{1}{2\varepsilon_K} \laplace_{{\mathbf y}^K} \left( G(\vartheta)^2 f^1 \right),
\]
where we denoted
\(  \label{def:Ff1}
    F[f^1]({\mathbf x}, {\mathbf y}^1) := 
      \int_{\R^{d(K+1)+1}} \grad W({\mathbf x} - {\mathbf x}_\ast) \cdot \left({\mathbf y}^1 - {\mathbf y}^1_\ast \right)
        f^1_\ast       \,\d{\mathbf x}_\ast \d {\mathbf y}_\ast \d \vartheta_\ast,
\)
with $f^1_\ast := f^1(t,{\mathbf x}_\ast, {\mathbf y}_\ast, \vartheta_\ast)$.
The weak formulation of the Vlasov equation with test function $\psi=\psi({\mathbf x}, {\mathbf y}, \vartheta)$ reads
\[
   \tot{}{t} \int f^1 \psi = \int \left( {\mathbf y}^1 f^1 \right) \cdot\grad_{\mathbf x} \psi
   +  \sum_{k=1}^{K-1}\frac{1}{\varepsilon_k} \int  \left( ({\mathbf y}^{k+1} - \alpha_k {\mathbf y}^{k}) f^1 \right) \cdot \grad_{{\mathbf y}^k} \psi \\
   - \frac{1}{\varepsilon_K} \int \left( \alpha_K {\mathbf y}^{K} f^1 \right) \cdot \grad_{{\mathbf y}^K} \psi \\
      + \int \left( F[f^1] f^1 \right) \part{\psi}{\vartheta} + \frac{1}{2\varepsilon_K} \int \left( G(\vartheta)^2 f^1 \right) \laplace_{{\mathbf y}^K} \psi,
\]
where all integrals are with respect to the variables $({\mathbf x}, {\mathbf y}, \vartheta)\in \R^{d(K+1)+1}$.
We now introduce the ansatz
\(    \label{ansatz}
    f^1(t, {\mathbf x}, {\mathbf y}, \vartheta) = f(t, {\mathbf x}, {\mathbf y})
       \,\delta(\vartheta - W\ast\varrho)
\)
with $\delta$ the Dirac delta function,
\[
   W\ast\varrho(t,{\mathbf x}) := \int_{\R^d} W({\mathbf x} - {\mathbf x}_\ast) \varrho(t, {\mathbf x}_\ast) \,\d {\mathbf x}_\ast
\]
and the particle number density
\begin{equation}   
\varrho(t,{\mathbf x}) := \int_{\R^{dK+1}}
 f^1(t, {\mathbf x}, {\mathbf y},\vartheta)  
\, \d {\mathbf y} \d\vartheta\,.
\label{def:rho}
\end{equation}
Inserting \eqref{ansatz} into the  weak formulation of the Vlasov equation and denoting
\[
   \varphi(t,{\mathbf x}, {\mathbf y}) := \psi({\mathbf x}, {\mathbf y}, W\ast\varrho(t,{\mathbf x}))
\]
gives
\(
   \tot{}{t} \int_\bullet f \varphi = \int \left( {\mathbf y}^1 f^1 \right) \cdot\grad_{\mathbf x} \psi
   +  \sum_{k=1}^{K-1} \frac{1}{\varepsilon_k} \int_\bullet  \left( ({\mathbf y}^{k+1} - \alpha_k {\mathbf y}^{k}) f \right) \cdot \grad_{{\mathbf y}^k} \varphi  \nonumber\\
   - \frac{1}{\varepsilon_K} \int_\bullet \left( \alpha_K {\mathbf y}^{K} f \right) \cdot \grad_{{\mathbf y}^K} \varphi   \label{eq:vlasov}\\
      + \int \left( F[f^1] f^1 \right) \part{\psi}{\vartheta} +  \frac{1}{2\varepsilon_K} \int_\bullet \left( G(W\ast\varrho)^2 f \right) \laplace_{{\mathbf y}^K} \varphi,  \nonumber
\)
where here and in the sequel the symbol $\int_\bullet$ denotes integration with respect to the variables $({\mathbf x}, {\mathbf y})\in \R^{d(K+1)}$,
while the plain integral symbol means integration with respect to $({\mathbf x}, {\mathbf y}, \vartheta)\in \R^{d(K+1)+1}$.
With the chain rule for the derivative we have
\[
   \grad_{\mathbf x} \varphi = \grad_{\mathbf x} \psi + \grad_{\mathbf x} (W\ast\varrho) \part{\psi}{\vartheta}.
\]
Consequently, the first term on the left-hand side of \eqref{eq:vlasov} reads
\(    \label{eq:tc}
   \int \left( {\mathbf y}^1 f^1 \right) \cdot\grad_{\mathbf x} \psi = \int_\bullet \left( {\mathbf y}^1 f \right) \cdot\grad_{\mathbf x} \varphi
      - \int \left( {\mathbf y}^1 f^1 \right) \cdot \grad_{\mathbf x} (W\ast\varrho) \part{\psi}{\vartheta}.
\)
Moreover, using \eqref{def:Ff1}, we write
\[
   \int \left( F[f^1] f^1 \right) \part{\psi}{\vartheta}  = 
      \iint \grad W({\mathbf x} - {\mathbf x}_\ast) \cdot \mathbf{y}^1 f^1_\ast f^1 \part{\psi}{\vartheta} -
      \iint \grad W({\mathbf x} - {\mathbf x}_\ast) \cdot \mathbf{y}^1_\ast f^1_\ast f^1 \part{\psi}{\vartheta},
\]
where the double integral sign $\iint$ means integral with respect to the variables
$({\mathbf x}, {\mathbf y}, \vartheta,{\mathbf x}_\ast, {\mathbf y}_\ast, \vartheta_\ast)\in \R^{2d(K+1)+2}$.
We now observe that the first term on the right-hand side
\[
   \iint \grad W({\mathbf x} - {\mathbf x}_\ast) \cdot \mathbf{y}^1 f^1_\ast f^1 \part{\psi}{\vartheta} 
   =
   \int \grad_{\mathbf x} (W\ast\rho) \cdot \mathbf{y}^1 f^1 \part{\psi}{\vartheta}    
\]
cancels the last term in \eqref{eq:tc}.
Moreover, the co-ordinate change $({\mathbf x}, {\mathbf y}, \vartheta) \rightleftarrows
({\mathbf x}_\ast, {\mathbf y}_\ast, \vartheta_\ast)$, using
the antisymmetry $\grad W({\mathbf x} - {\mathbf x}_\ast) = - \grad W({\mathbf x}_\ast - {\mathbf x})$, gives
\[
   - \iint \grad W({\mathbf x} - {\mathbf x}_\ast) \cdot \mathbf{y}^1_\ast f^1_\ast f^1 \part{\psi}{\vartheta}
   =
   \iint \grad W({\mathbf x} - {\mathbf x}_\ast) \cdot \mathbf{y}^1 f^1 f^1_\ast \part{\psi_\ast}{\vartheta}  \\
   =
   \int_\bullet \grad_{\mathbf x} \left[ W\circledast \left(f \part{\psi}{\vartheta} \right) \right] \cdot \mathbf{y}^1 f,
\]
where we denoted
\[
   W\circledast \left(f \part{\psi}{\vartheta} \right)(t,{\mathbf x}) :=
      \int_{\R^{d(K+1)}} W({\mathbf x} - {\mathbf x}_\ast) f(t, {\mathbf x}_\ast, {\mathbf y})
         \part{\psi}{\vartheta}({\mathbf x}_\ast, {\mathbf y}, W\ast\varrho(t,{\mathbf x}))
          \, \d {\mathbf y} \d {\mathbf x}_\ast.
\]
We conclude that
\[
   \int \left( {\mathbf y}^1 f^1 \right) \cdot\grad_{\mathbf x} \psi + \int \left( F[f^1] f^1 \right) \part{\psi}{\vartheta}
   =
   \int_\bullet \left( {\mathbf y}^1 f \right) \cdot\grad_{\mathbf x} \left[ \varphi +  W\circledast \left(f \part{\psi}{\vartheta} \right) \right].
\]
Finally, for the time derivative term in \eqref{eq:vlasov} we again use the chain rule
\[
    \part{\varphi}{t} = \part{\psi}{\vartheta} \left( W\ast\part{\varrho}{t} \right)
\]
and obtain
\[
   \tot{}{t} \int_\bullet f \varphi &=& 
   \int_\bullet \part{f}{t} \varphi + \int_\bullet f \part{\psi}{\vartheta} \left( W\ast\part{\varrho}{t} \right) \\
   &=&
    \int_\bullet \part{f}{t} \varphi + \int_\bullet \part{\varrho}{t} \, W\ast \left( f \part{\psi}{\vartheta} \right) \\
   &=&
    \int_\bullet \part{f}{t} \left[ \varphi +  W\circledast \left(f \part{\psi}{\vartheta} \right) \right].
\]
Therefore, denoting $\xi:= \varphi +  W\circledast \left(f \part{\psi}{\vartheta} \right)$ and noting that
$\grad_{\mathbf{y}} \xi = \grad_{\mathbf{y}} \varphi$, we write \eqref{eq:vlasov} as
\[
   \int_\bullet \part{f}{t} \xi = \int_\bullet \left( {\mathbf y}^1 f \right) \cdot\grad_{\mathbf x} \xi
      +  \sum_{k=1}^{K-1} \frac{1}{\varepsilon_k} \int_\bullet  \left( ({\mathbf y}^{k+1} - \alpha_k {\mathbf y}^{k}) f \right) \cdot \grad_{{\mathbf y}^k} \xi  \nonumber\\
   - \frac{1}{\varepsilon_K} \int_\bullet \left( \alpha_K {\mathbf y}^{K} f \right) \cdot \grad_{{\mathbf y}^K} \xi
   + \frac{1}{2\varepsilon_K} \int_\bullet \left( G(W\ast\varrho)^2 f \right) \laplace_{{\mathbf y}^K} \xi.  \nonumber
\]
This is the weak formulation of the Fokker-Planck equation
\begin{eqnarray}
\partial_t f \!\!&+&\!\! {\mathbf y}^1\cdot \grad_{\mathbf x} f 
\,+\, \sum_{k=1}^{K-1} \frac{1}{\varepsilon_k} \grad_{{\mathbf y}^k}\cdot \left( ({\mathbf y}^{k+1} - \alpha_k {\mathbf y}^k) f \right) \nonumber \\
&=& \frac{1}{\varepsilon_K} \, \grad_{{\mathbf y}^K}\cdot \left( \alpha_k {\mathbf y}^K f + \frac{G(W\ast\varrho)^2}{2} \grad_{{\mathbf y}^K}  f \right),
\label{eq:FP}
\end{eqnarray}
with the particle number density $\varrho=\varrho(t,{\mathbf x})$ given by \eqref{def:rho}, i.e.,
\begin{equation}   
\varrho(t,{\mathbf x}) = \int_{\R^{dK}}
 f(t, {\mathbf x}, {\mathbf y})  
\, \d {\mathbf y}.
\label{def:FPrho}
\end{equation}


\section{Equilibrium patterns in the large population limit}
\label{sec:patterns}

Equilibrium patterns in the limit $N\to\infty$ of the discrete system~\eqref{SDE1} 
correspond to steady states of the Fokker-Planck equation \eqref{eq:FP}--\eqref{def:FPrho}.
In Sections~\ref{sec5} and~\ref{sec6} we will present results of stochastic simulations of the discrete system
where the positions of individuals, ${\mathbf x}_i = {\mathbf x}_i(t) \in \R^d$, $i\in[N]$, evolve in the $d$-dimensional cube
\begin{equation}
\Omega = (0,1)^d \qquad \mbox{with periodic boundary conditions.}
\label{domainomega}
\end{equation}
We therefore consider equation~\eqref{eq:FP} posed on a torus $\Omega$ in the ${\mathbf x}$ variable, 
while the ${\mathbf y}^k$ variables, for $k\in[K]$, take their values in the full space $\R^d$.
Since for any fixed ${\mathbf x} \in\Omega$ the function 
$$
{\mathbf y}^K \mapsto \exp\left( - \frac{\alpha_K |{\mathbf y}^K|^2}{G(W\ast\varrho({\mathbf x}))^2} \right)
$$
is an equilibrium for the Fokker-Planck operator on the right-hand side of equation~\eqref{eq:FP},
we have the following homogeneous ({i.e}., constant in ${\mathbf x}$) steady state solution of the
Fokker-Planck equation~\eqref{eq:FP}
\begin{equation}   
f(t, {\mathbf x}, {\mathbf y}) 
= C \, 
\exp\left( - \frac{\alpha_K |{\mathbf y}^K|^2}{G(\bar\varrho)^2} \right) \prod_{k=1}^{K-1} 
      \delta \!\left({\mathbf y}^k - \alpha_k^{-1} {\mathbf y}^{k+1} \right),
\label{eq:steady}
\end{equation}
where $\bar\varrho>0$ is the constant particle density in $\Omega$, $\delta$ is the Dirac delta function on $\R^d$ 
and $C>0$ is the normalization constant. Using~(\ref{def:FPrho}) and~\eqref{eq:steady}, we get
\begin{equation*}   
\bar\varrho  = \int_{\R^{dK}} 
 f(t, {\mathbf x}, {\mathbf y})  
\, \d {\mathbf y}
     \,=\,  C\, \left( \!\frac{\pi}{\alpha_K} \!\right)^{\!\!d/2} G(\bar\varrho)^d\,.
\end{equation*}
Consequently, the normalization constant $C$ in equation~\eqref{eq:steady} is given by
$$
C \,=\, \bigg( \!\frac{\alpha_K}{\pi} \! \bigg)^{\!\!d/2} \frac{\bar\varrho}{G(\bar\varrho)^d} \,.
$$
In particular, equation~\eqref{eq:steady} represents a steady state that is constant in ${\mathbf x}$, concentrated 
on the diagonals ${\mathbf y}^k \,=\, \alpha_k^{-1}\, {\mathbf y}^{k+1}$, for $k\in[K-1]$, and Gaussian distributed 
in the ${\mathbf y}^K$~variable. From the point of view of the discrete individual-based model~\eqref{SDE1},
this corresponds to the equilibrium where no aggregation takes place (constant particle density) and the internal 
variables ${\mathbf y}^k$, for $k\in[K]$, are subject to the Brownian motion.
To characterise the steady states which are non-constant in the ${\mathbf x}$~variable,
we first write the momentum system corresponding to the Fokker-Planck equation~\eqref{eq:FP}.
Integrating~\eqref{eq:FP} with respect to all ${\mathbf y}^k$~variables, we obtain
\begin{equation}
\partial_t \varrho \,+\, \grad_{\mathbf x} \cdot \int_{\R^{dK}} {\mathbf y}^1 \, f \,\d {\mathbf y} \,=\, 0.
\label{eq:m0}
\end{equation}
We observe that constructing $f$ radially symmetric in ${\mathbf y}^1$ ({i.e.}, depending on ${\mathbf y}^1$ 
only through the modulus $|{\mathbf y}^1|$) annihilates the term $\int_{\R^{dK}} {\mathbf y}^1 \, f \,\d {\mathbf y}$
and equation~\eqref{eq:m0} reduces to $\partial_t \varrho = 0$.
Multiplying equation~\eqref{eq:FP} by ${\mathbf y}^K$ and integrating with respect to all ${\mathbf y}^k$ variables,
we obtain
\begin{equation}
\partial_t \left( \varrho \, m^K \right) 
\,+\,  
\grad_{\mathbf x}
\cdot \int_{\R^{dK}} \left( {\mathbf y}^1 \otimes {\mathbf y}^K \right) f  \,\d {\mathbf y} 
\, = \,
- \frac{1}{\varepsilon_K} \, \alpha_K \, \varrho \, m^K \, ,
\label{eq:mm}
\end{equation}
where we have denoted $\varrho \, m^K := \int_{\R^{dK}} {\mathbf y}^K \, f \,\d {\mathbf y}$.
From here we infer that in equilibrium we either have empty regions where $\varrho=0$,
or regions with positive density $\varrho>0$, but zero moment $m^K=0$.
Moreover, let us observe that any equilibrium of \eqref{eq:FP} needs to be concentrated
of the diagonals ${\mathbf y}^k = \alpha_k^{-1} \, {\mathbf y}^{k+1}$, for $k\in[K-1]$.
Since we have imposed radial symmetry in ${\mathbf y}^1$, this implies that the equilibrium 
is radially symmetric in all ${\mathbf y}^k$ variables. Then an easy calculation reveals that
$$
\left( \prod_{k=1}^{K-1} \alpha_k \!\right)
\grad_{\mathbf x} \cdot \int_{\R^{dK}} \left({\mathbf y}^1 \otimes {\mathbf y}^K \right) f  \,\d {\mathbf y}
\, = \,
\grad_{\mathbf x} \cdot  \int_{\R^{dK}} \left( {\mathbf y}^K \otimes {\mathbf y}^K \right) f  \,\d {\mathbf y}
\, = \,
\frac{2}{d} \, \grad_{\mathbf x}  \!\left(\varrho \, e^K \right),
$$
where we have denoted $\varrho \, e^K = \frac12 \int_{\R^{dK}} \big|{\mathbf y}^K\big|^2 f \, \d {\mathbf y}$.
We see that $\varrho \, e^K$ plays the role of pressure in the momentum equation~\eqref{eq:mm}
and in equilibrium we must have $\varrho \, e^K \equiv \mbox{constant}$ on $\Omega$.
Finally, multiplying equation~\eqref{eq:FP} by $\big|\mathbf y^K\big|^2$ and integrating by parts, we obtain
$$
\partial_t \left( \varrho \, e^K \right) 
\,+\, 
\grad_{\mathbf x} \cdot \left( \frac12 \int_{\R^{dK}} \big|{\mathbf y}^K\big|^2 {\mathbf y}^1 f \, \d {\mathbf y} \right)
\,=\, 
- \frac{1}{\varepsilon_K} \, \alpha_K \, \varrho \, e^K \, + \, \frac{d}{2\varepsilon_K} \, G(W\ast\varrho)^2 \, \varrho.
$$
We observe that in equilibrium we must have
$$
\frac{d}{2} \, G(W\ast\varrho)^2 \, \varrho \, = \, \alpha_K \, \varrho \, e^K \,.
$$
Due to the requirement that $\varrho \, e^K$ be constant in $\Omega$, we finally arrive
at the condition characterising the equilibrium density profiles,
namely,
\begin{equation}
G(W\ast\varrho)^2\, \varrho \, \equiv \, C_0 \qquad \mbox{in } \Omega\, ,
\label{eq:eqrho}
\end{equation}
where the constant $C_0>0$ is determined by the initial mass.
This allows for nonhomogeneous number density profiles~$\varrho=\varrho({\mathbf x})$,
i.e., presence of clusters.%
\begin{figure}
\begin{center}
\includegraphics[width=.45\textwidth]{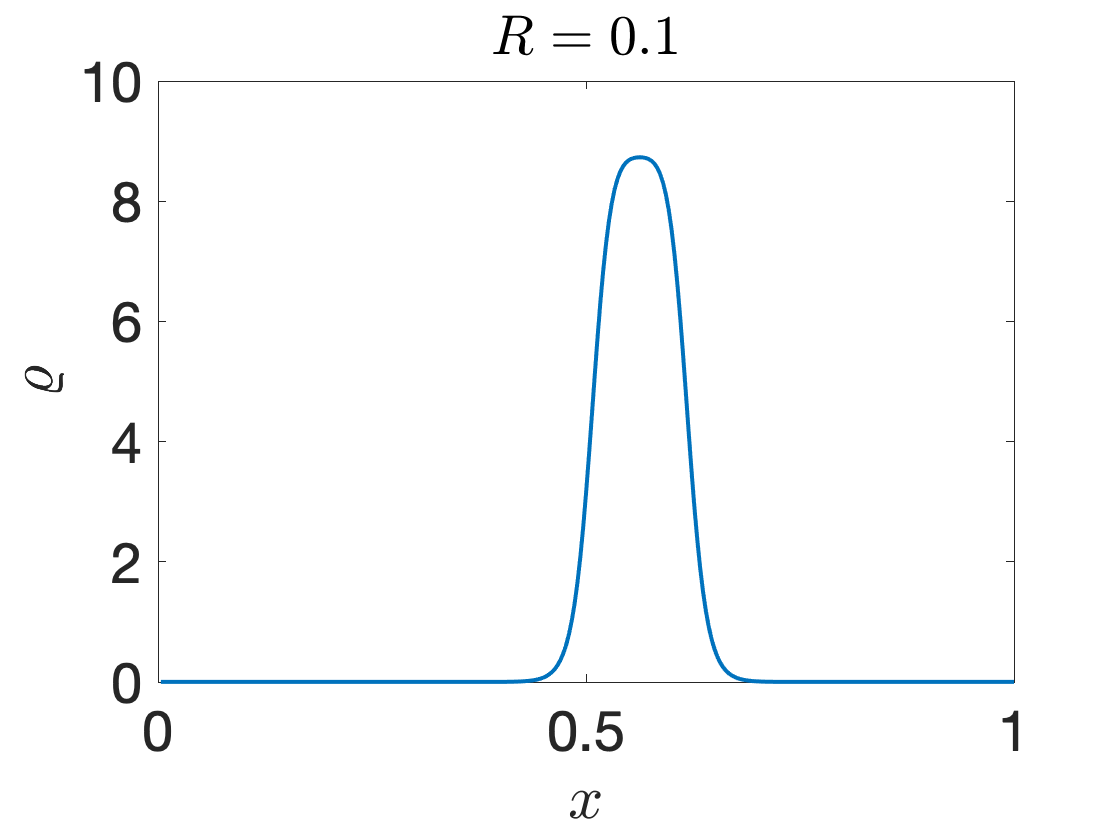}
\includegraphics[width=.45\textwidth]{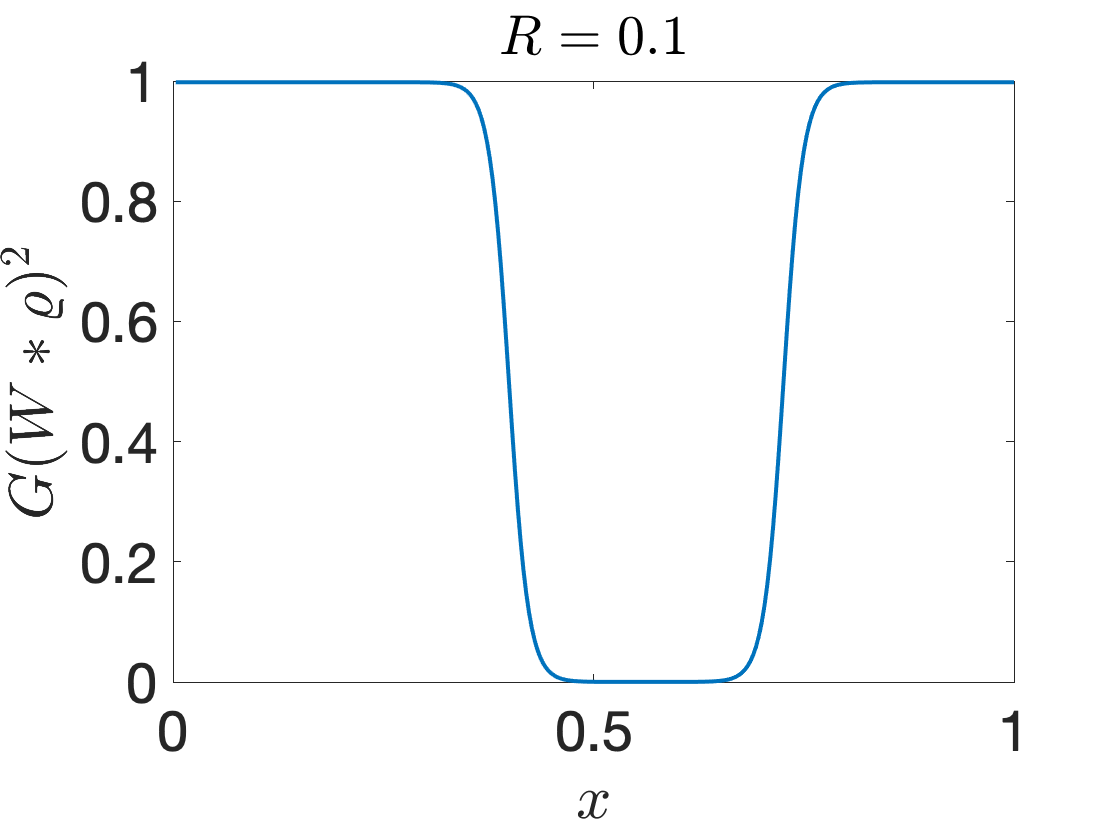}\\
\includegraphics[width=.45\textwidth]{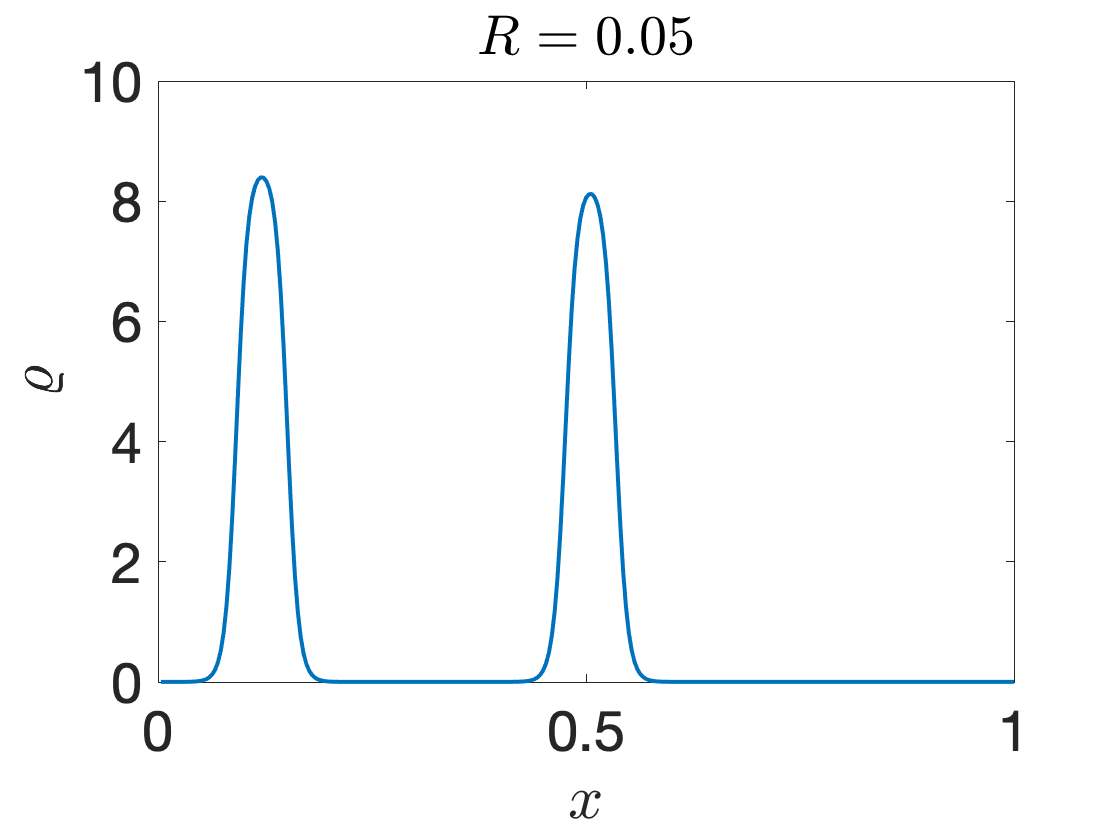}
\includegraphics[width=.45\textwidth]{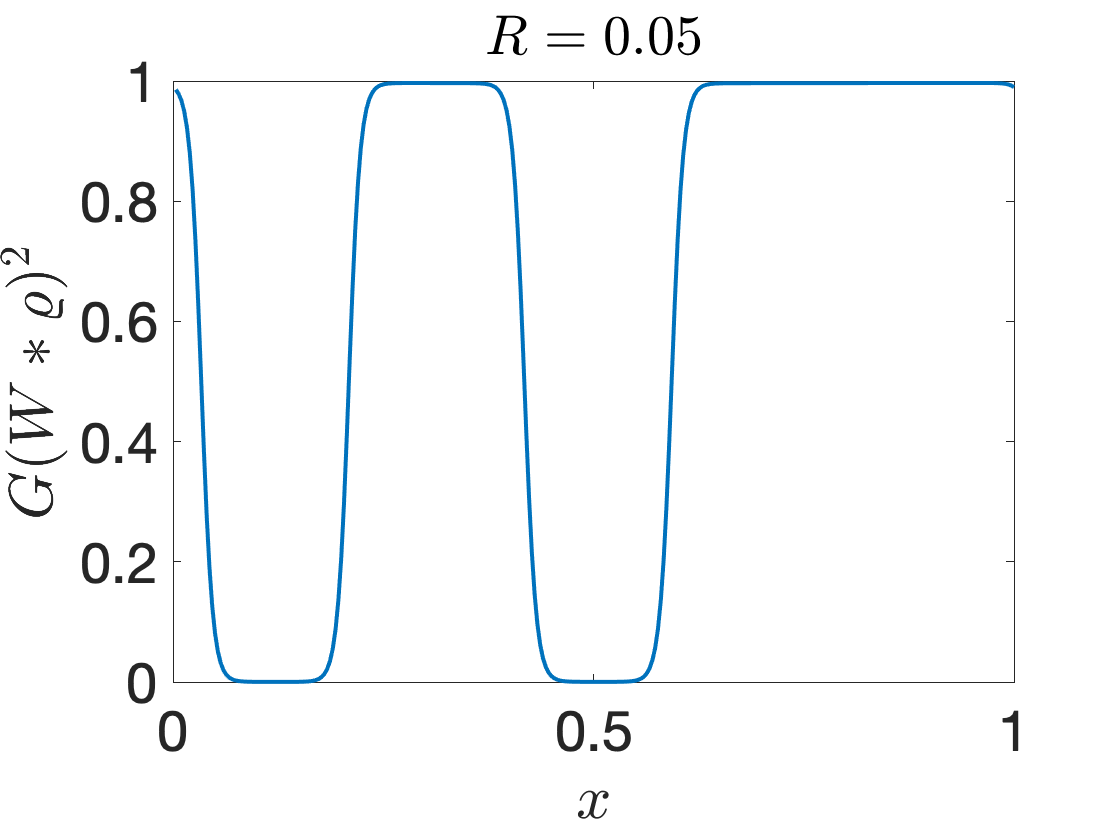}\\
\includegraphics[width=.45\textwidth]{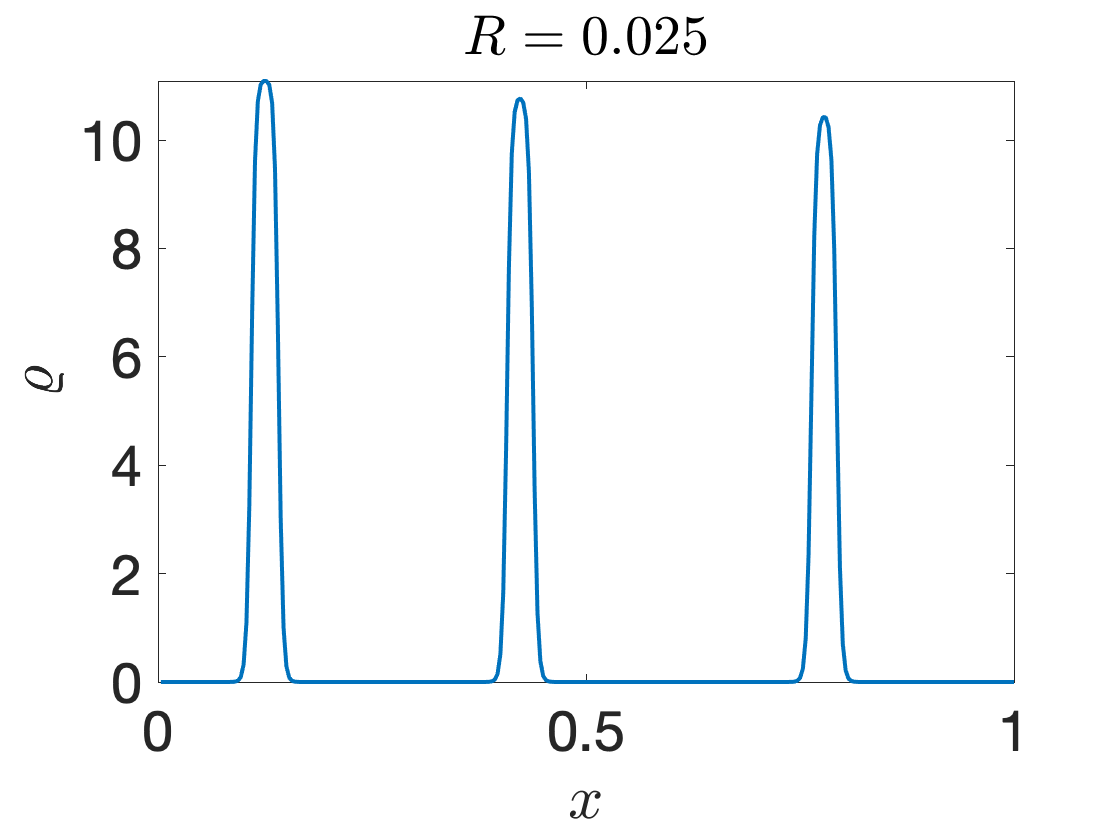}
\includegraphics[width=.45\textwidth]{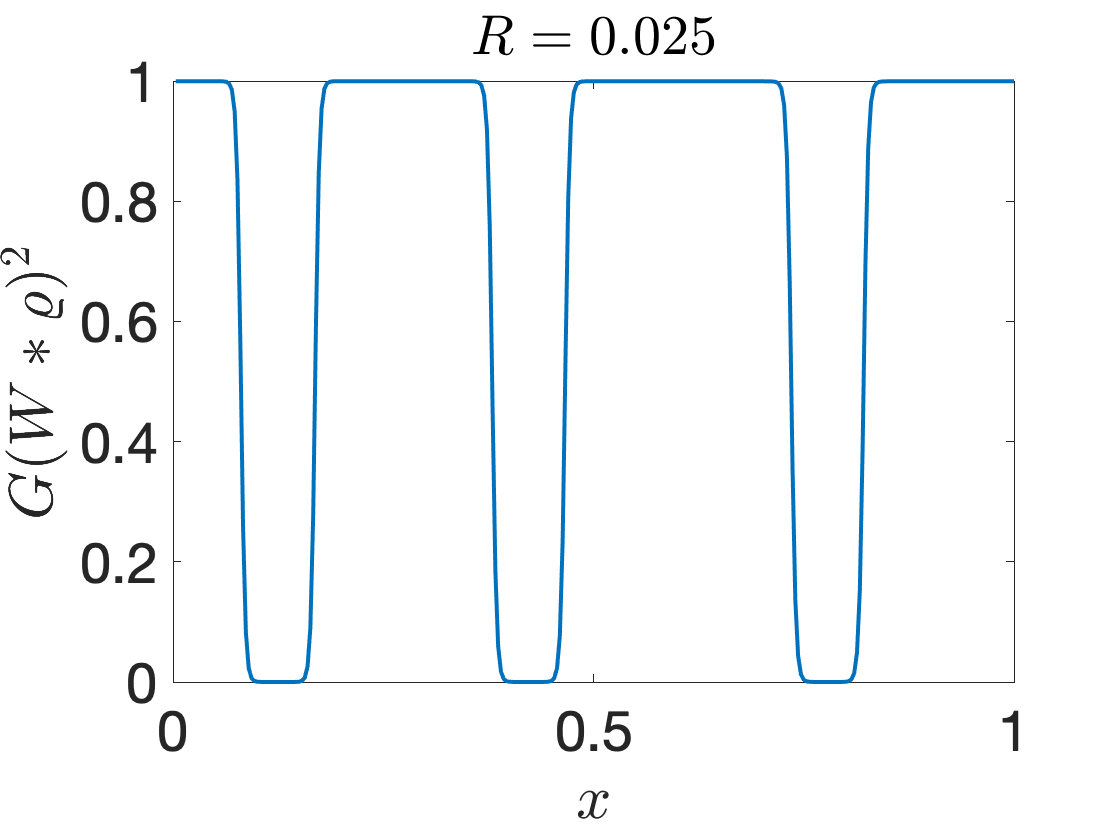}
\end{center}
\caption{{\it One-dimensional equilibrium profiles satisfying the condition $(\ref{eq:eqrho})$,
obtained by solving equation~$\eqref{eq:rho}$ in domain~$(\ref{domainomega})$ for $d=1$,
subject to the initial datum $\varrho(t=0,x) = 1+ 10^{-1} \sin(2\pi x(2-x))^2$,
until a steady state is reached. We use $G(s)$ and $W(x)$ given by~$(\ref{choiceGandW})$
with $R=0.1$ $($top$)$, $R=0.05$ $($middle$)$ and $R=0.025$ $($bottom$)$.
The left panels depict the steady state density $\varrho=\varrho(x)$, while the right panels 
visualize the function $G(W\ast\varrho)^2$. 
}}
\label{figure1}
\end{figure}
To numerically compute such nonhomogeneous profiles~$\varrho=\varrho({\mathbf x})$ satisfying condition~\eqref{eq:eqrho},  
we can solve equation~\eqref{eq:rho} in the unstable regime, subject to
an initial datum that is a perturbation of the constant steady state.
We employ this strategy in the spatially one-dimensional setting, {i.e.}, for $d=1$ in our
domain~(\ref{domainomega}). The results are presented in Figure~\ref{figure1}, where we use 
\begin{equation}
G(s) := e^{-s} 
\qquad \mbox{and} \qquad 
W(x) := \frac{\chi_{[0,R]}(|x|)}{2R} \quad \mbox{with some} \;\;R>0\,,
\label{choiceGandW}
\end{equation}
where $\chi_{[0,R]}$ denotes the characteristic function of the interval $[0,R]$,
{i.e.}, the kernel $W=W(x)$ corresponds to the sampling radius $R>0$.
Using the semi-implicit finite difference discretization for the space variable
and first-order forward Euler method for the time variable, 
the steady states are plotted for $R\in \{0.1, 0.05, 0.025\}$ in Figure~\ref{figure1}.
Let us note that the steady states are not unique and the particular pattern
produced by the simulation depends on the choice of the initial datum;
for all three simulations presented in Figure~\ref{figure1} we used $\varrho(t=0,x) = 1+ 10^{-1}\sin(2\pi x(2-x))^2$.
The same approach can be applied in higher-dimensional settings as well,
see~\cite[Section 5.4]{BHW:2012:PhysD} for the results with $d=2$.

Finally, let us mention that in the formal limit as $R \to 0+$ in \eqref{choiceGandW}
we obtain $W(x) = \delta(x)$, the Dirac delta function in $\R^d$. Then $W\ast\varrho = \varrho$
and equation~\eqref{eq:rho} becomes a nonlinear but local reaction-diffusion equation.
However, it does not produce any nontrivial patterns, since the following alternative holds:
either $G=G(s)$ is such that all homogeneous steady states are asymptotically stable,
or equation~\eqref{eq:rho} is ill-posed. We refer to~\cite[Section 4.4]{BHW:2012:PhysD} for more details.

\section{Spontaneous aggregation in one spatial dimension ($d=1$)}
\label{sec5}

\noindent
As the dimensionality of the Fokker-Planck equation~\eqref{eq:FP} becomes prohibitive for numerical simulations
even with moderate values of $K$, we use stochastic simulations of the individual-based model
given by equations~\eqref{rho_i} and~\eqref{SDE1} to systematically investigate the impact of the number 
of memory layers $K\geq 1$ on the clustering properties of the spontaneous aggregation model.
We use $\varepsilon_k$ and $\alpha_k$, for $k\in[K]$, given by equations~(\ref{epsare1})
and~(\ref{alphakarealpha}), respectively, where we choose $\alpha=1$ in equation~(\ref{alphakarealpha}).
The distance between agents is calculated over the torus, {i.e.}, taking into account the periodic boundary 
conditions on $\Omega$. We simulate $N=400$ agents moving in the domain~(\ref{domainomega}) for $d=1$. 
The response function $G(s)$ and the interaction kernel $W(x)$ are given by~(\ref{choiceGandW})
with the sampling radius $R=1/40=0.025$. 

We initialize the simulation by randomly generated agent positions $x_i\in\Omega$, for $i\in[N]$,
using uniformly distributed initial positions in $\Omega$. All internal variables are initialized as zeros,
{i.e.}, $y_i^k(0)=0$ for all $i\in[N]$ and $k\in[K]$. We discretize equations~\eqref{SDE1} using 
the Euler-Maruyama scheme with timestep $\Delta t=10^{-3}$. We calculate $100$ stochastic realizations 
of the individual-based model in the time interval $[0,10^3]$, {i.e.}, we calculate the time evolution 
over $10^6$ timesteps for each realization.
\rev{
Initially, clusters emerge from random fluctuations in the particle distribution.
Subsequently, small clusters dissociate and a dynamic equilibrium is established between clustered particles and freely moving outliers.
This equilibrium typically forms on a relatively fast timescale, e.g.,
at the order of the first $10^5$ timesteps (out of the total $10^6$ timesteps carried out in the simulation).
Then, the number of clusters remains constant until the end of the simulation.
The location, shape, and size of the clusters fluctuate only minimally.
We thus observe the emergence of well-defined quasi-stationary states. 
}
We record the particle positions at the final time $t=10^3$ 
and use these for evaluating clustering properties of the model.

\begin{figure}[t]
\begin{center}
\includegraphics[width=.33\textwidth]{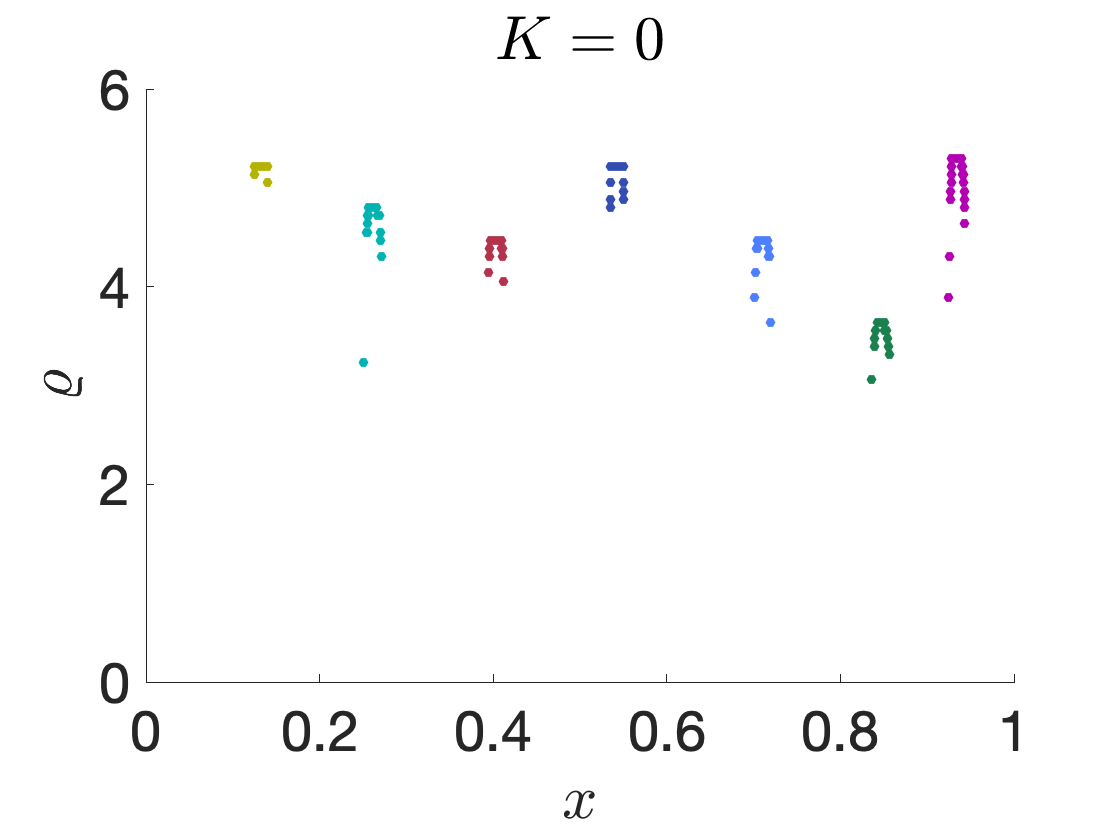}%
\includegraphics[width=.33\textwidth]{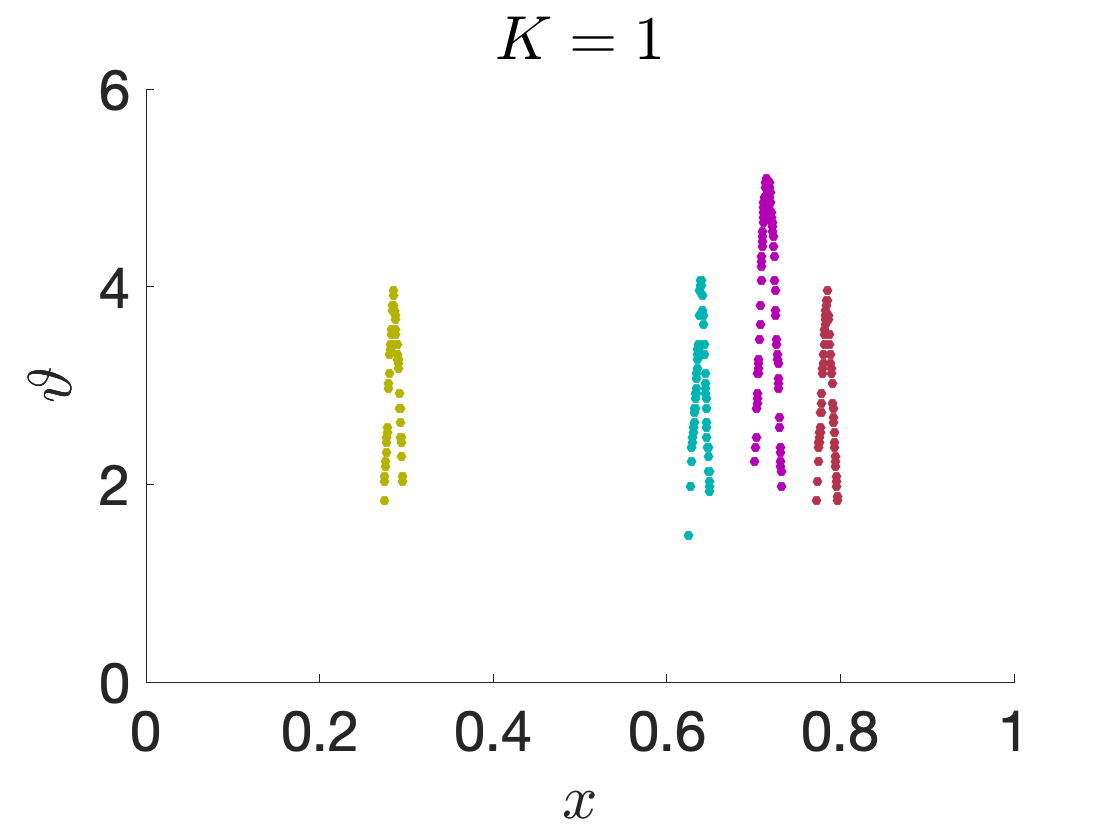}%
\includegraphics[width=.33\textwidth]{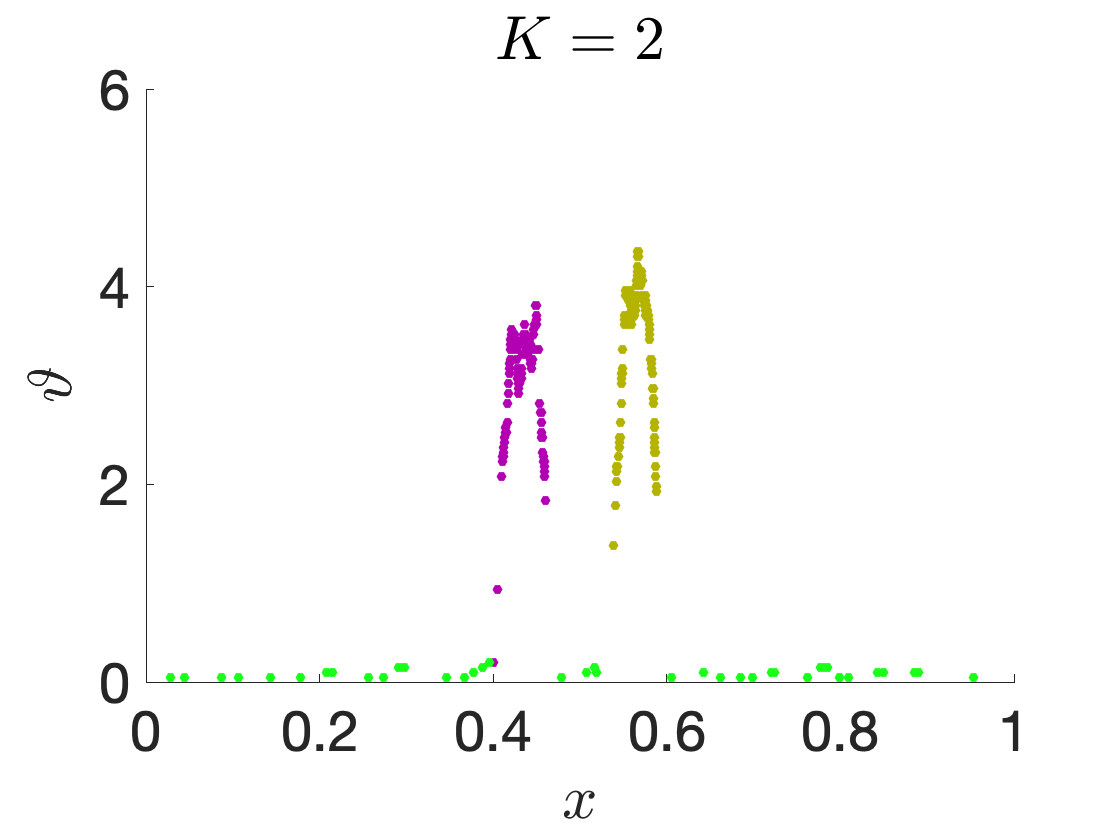}\\
\includegraphics[width=.33\textwidth]{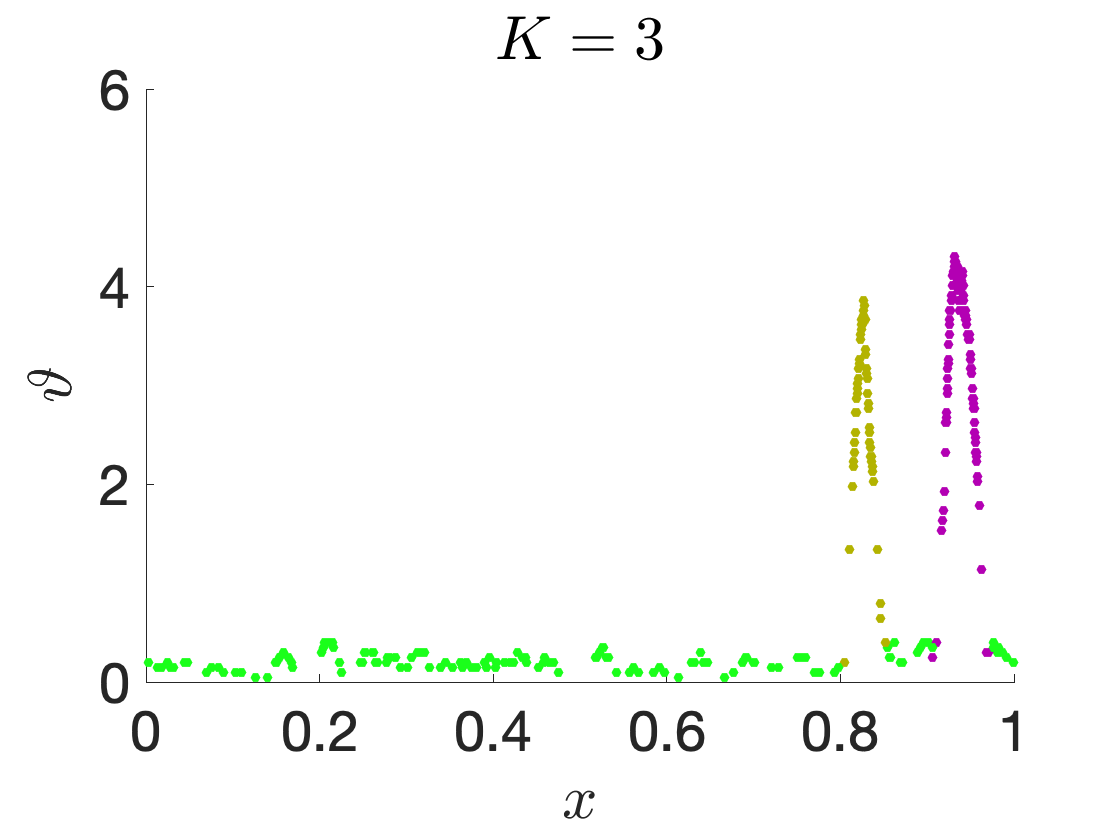}%
\includegraphics[width=.33\textwidth]{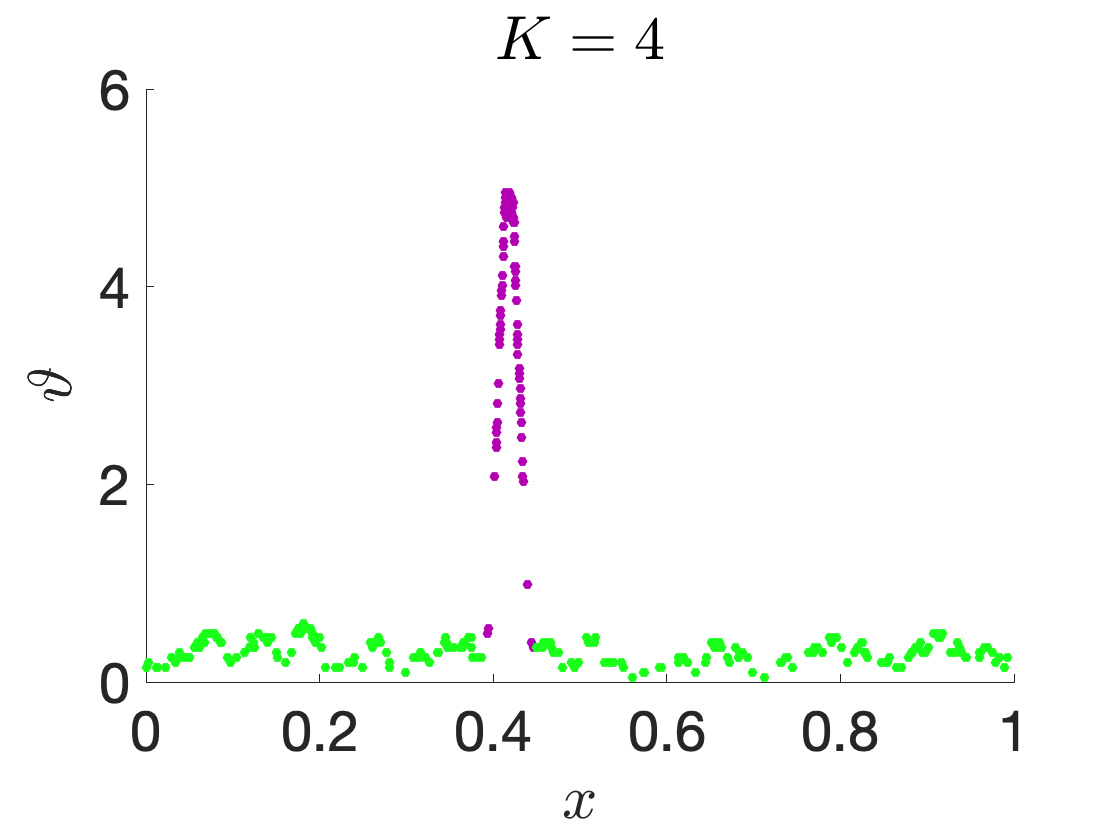}%
\includegraphics[width=.33\textwidth]{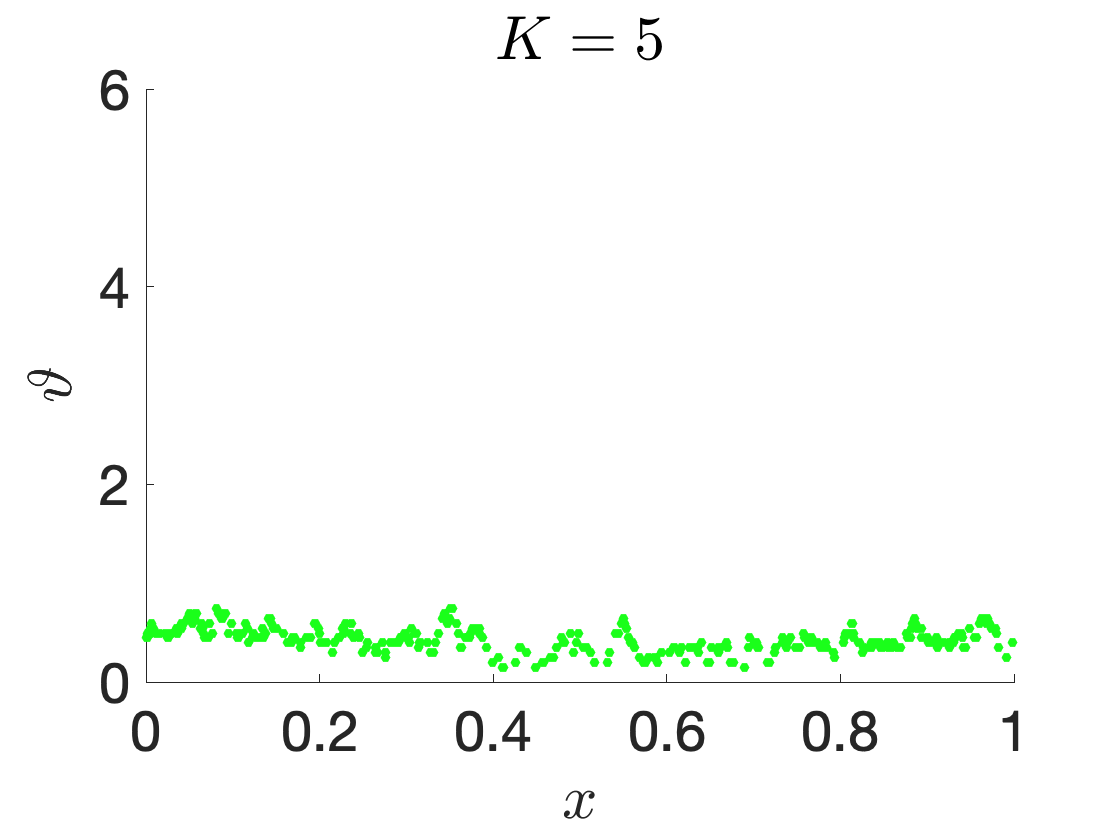}
\end{center}
\vskip -5mm
\caption{{\it Simulations of the individual-based model~$\eqref{rho_i}$ and~$\eqref{SDE1}$ for
$K\in\{1,\dots,5\}$, and the model~\eqref{rho_i}--\eqref{model1} without memory ($K=0)$.
We used $N=400$ agents moving in the domain $\Omega$ given by~$(\ref{domainomega})$ with $d=1$.
$G(s)$ and $W({\mathbf x})$ are given by~$(\ref{choiceGandW})$ with $R=0.025$.
For each value of $K\in\{0,1,\dots,5\}$, the plots capture the particle positions at the final time 
$t=10^3$ $($horizontal axis$)$ and the perceived density of their neighbours $\vartheta_i$ given by~\eqref{rho_i} 
$($vertical axis$)$. The clusters, differentiated by colour, are identified using the {\rm DBSCAN} method 
with parameters $\mathtt{epsilon}=0.025$ and $\mathrm{MinPoints}=20$. The light green points are outliers, 
{i.e.}, particles not belonging to any cluster.
}}
\label{figure2}
\end{figure}

For identification of the clusters we use the {\tt MATLAB} implementation
of the Density-based spatial clustering (DBSCAN) method~\cite{DBSCAN}.
As we expect the clusters to be of size comparable to the sampling radius
$R=0.025$, we set the parameter $\mathtt{epsilon}$, specifying the radius of a neighbourhood
with respect to some point of the DBSCAN method, to $\mathtt{epsilon}=0.025$.
Moreover, based on experimentation, we found that the choice of the second
parameter of DBSCAN, $\mathtt{minPts}=20$, leads to the best results
in identification of clusters. Examples of the results, recorded at the final 
timestep of the simulations with $K\in\{1,2,3,4,5\}$, are plotted 
in Figure~\ref{figure2}. For comparison, we also simulate the 
system~\eqref{rho_i}--\eqref{model1}, {i.e.}, spontaneous aggregation without memory; 
we refer to the corresponding results by $K=0$.

\begin{figure}
(a) \hskip 7.5cm (b) \hfill\break
\includegraphics[width=.49\textwidth]{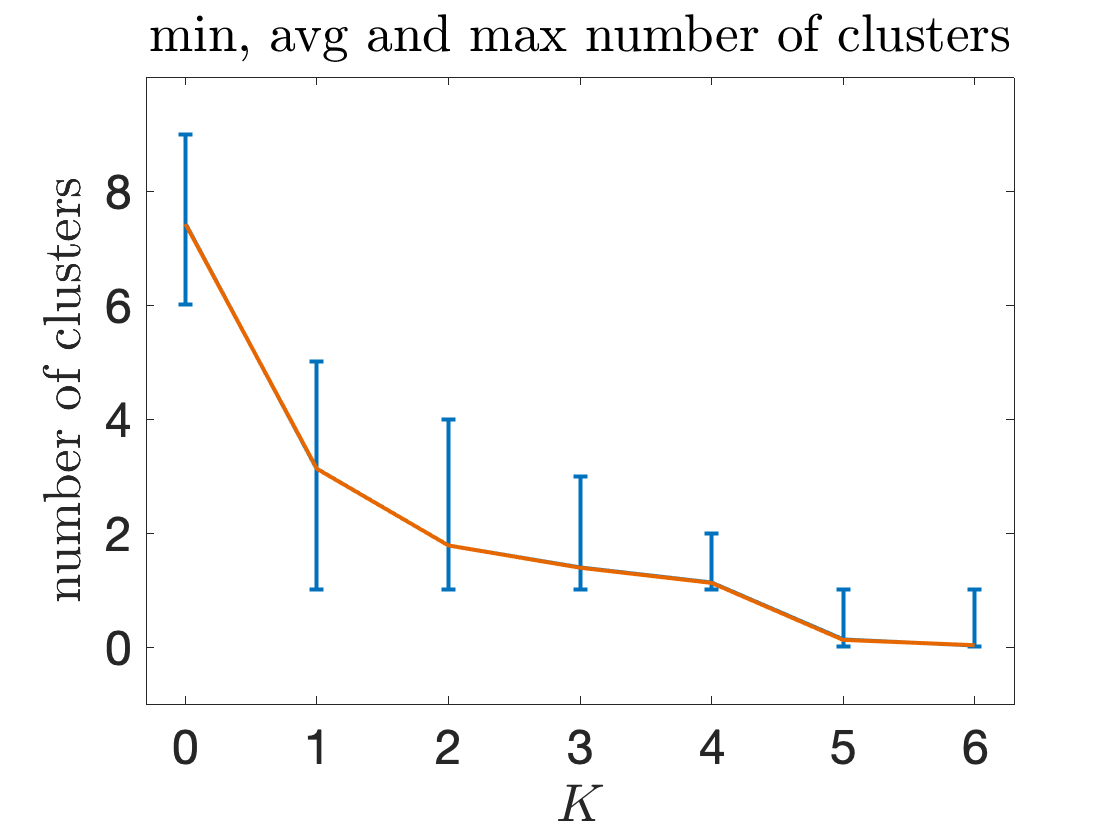}
 \includegraphics[width=.49\textwidth]{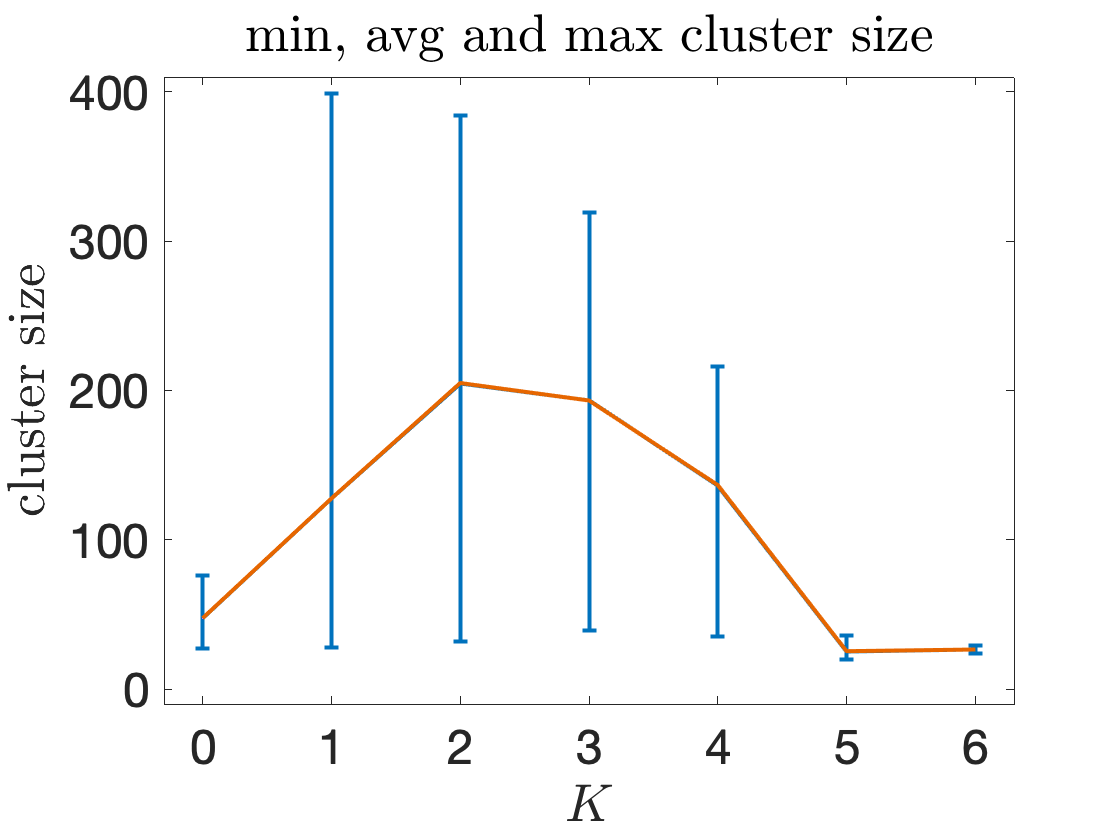} \hfill\break
(c) \hskip 7.5cm (d) \hfill\break 
\includegraphics[width=.49\textwidth]{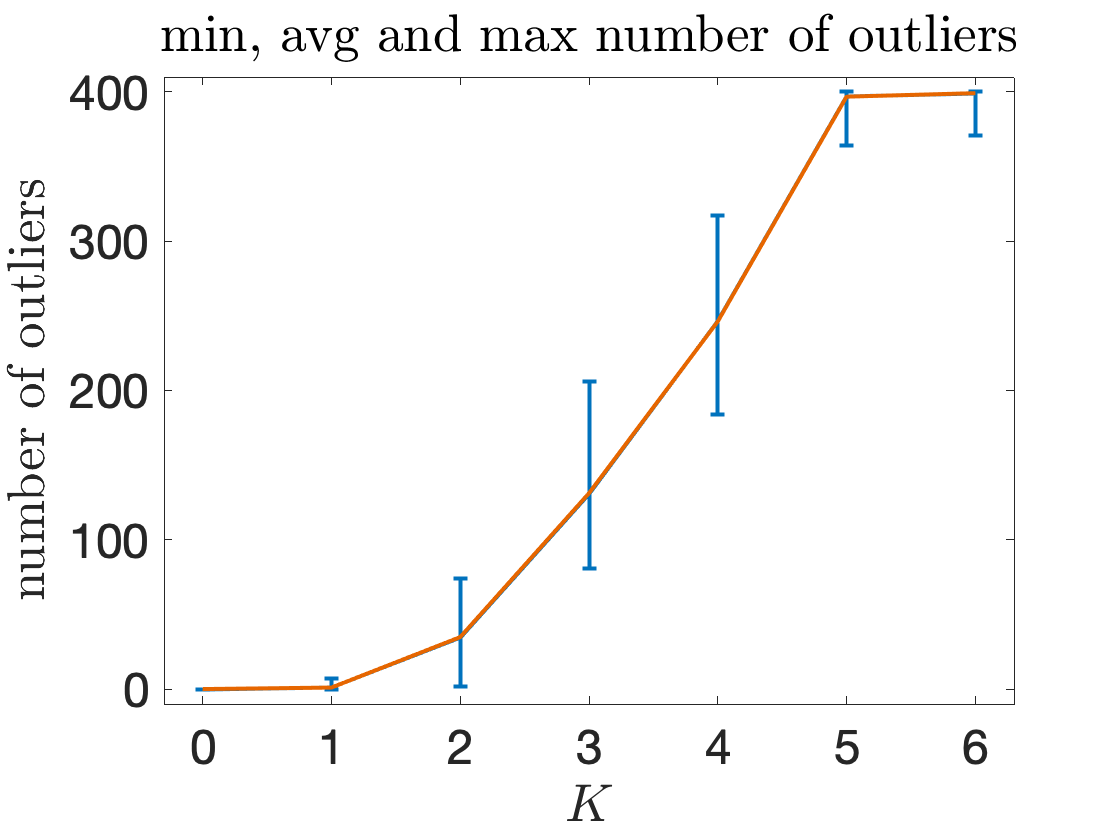}
\includegraphics[width=.49\textwidth]{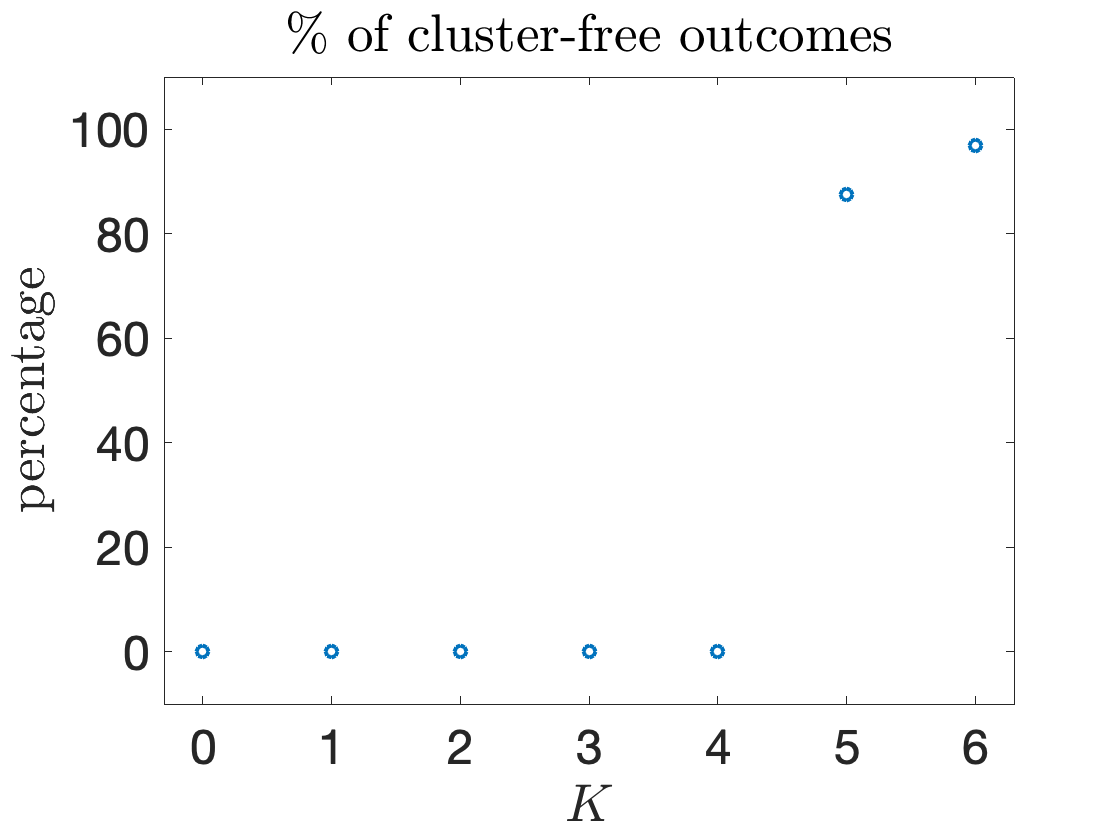}
\caption{{\it
Statistics of the clustering behaviour over $100$ realizations of the individual-based model given 
by equations~$\eqref{rho_i}$ and~$\eqref{SDE1}$ for dimension $d=1$, $N=400$ agents, $G(s)$ and 
$W(x)$ given by~$(\ref{choiceGandW})$ with $R=0.025$ and $K \in \{1, 2, \dots, 6\}$.
The case $K=0$ (no memory) refers to stochastic simulation of the system
\eqref{rho_i}--\eqref{model1}. Other parameters are the same as in Figure~$\ref{figure2}$. \hfill\break
} (a) {\it average $($orange solid line$)$, minimum and maximum $($indicated by blue error bars$)$ number of clusters identified in the 
final timestep of the simulations at time $t=10^3$, \hfill\break
} (b) {\it average $($orange solid line$)$, minimum and maximum $($blue error bars$)$ cluster sizes, \hfill\break
} (c) {\it average $($orange solid line$)$, minimum and maximum $($blue error bars$)$ number of outliers, {i.e.},
particles that do not belong to any cluster}, \hfill\break
(d) {\it percentage of simulation outcomes $($out of the $100$ runs$)$
that did not produce any clusters.}
}
 \label{figure3}
\end{figure}

In Figure~\ref{figure3}, we provide statistics of the clustering behaviour over $100$ realizations
of the individual-based stochastic model given by equations~\eqref{rho_i} and~\eqref{SDE1}, performed 
for each value $K \in \{1, 2, \dots, 6\}$, together with the results obtained by the
system~\eqref{rho_i}--\eqref{model1} (referred as the $K=0$ case).
In Figure~\ref{figure3}(a), we plot the number of clusters identified in the final
timestep of the simulation by the DBSCAN method. The solid line represents the average over the $100$ 
realizations, while the blue `error bars' indicate the minimum and maximum values. Here, we clearly observe 
a tendency toward the formation of fewer clusters as $K$ increases. 
Figure~\ref{figure3}(b) shows the minimum, average and maximum cluster sizes observed over the course 
of $100$ stochastic realizations. For $K\in\{0,1,2\}$ the size of the clusters increases with~$K$, and 
decreases for $K>3$. Figure~\ref{figure3}(c) shows the average (solid line)
and minimum/maximum (error bars) number of outliers, {i.e.}, particles that do not belong to 
any cluster. The number of outliers increases from almost zero for $K\in \{0,1\}$,
up to almost all particles being outliers for $K\in \{5,6\}$;
{c.f.} the corresponding panels in Figure~\ref{figure2}. Finally, the percentages
of the outcomes (out of the $100$ stochastic simulations)
that did not produce any clusters are plotted in Figure~\ref{figure3}(d). We observe that for 
$K\in\{0,1,2,3,4\}$ the percentage is zero (clusters have always been formed).
For $K=5$ the percentage increases steeply and is almost $100\%$ for $K=6$.

In summary, increasing memory length disrupts spontaneous aggregation.
In the next section, we present an analogous set of simulations in a two-dimensional spatial setting.
There, we also provide statistical evidence that memory weakens the particles' responsiveness to environmental stimuli.
We argue that this reduction in responsiveness accounts for the decreased tendency to cluster observed 
at higher values of~$K$.

\section{Spontaneous aggregation in two spatial dimensions ($d=2$)}

\label{sec6}

\noindent
In this section, we present results of stochastic simulations of the individual-based model 
given by equations~$\eqref{rho_i}$ and~$\eqref{SDE1}$ for $d=2$ (spatially two-dimensional setting),
where ${\mathbf x}_i={\mathbf x}_i(t)\in\Omega$ with $\Omega$ given by~(\ref{domainomega})
for $d=2.$ The internal variables ${\mathbf y}_i^k = {\mathbf y}_i^k(t)$ evolve in the full space $\R^2$.
As in Section \ref{sec5}, we set $\varepsilon_k$ and $\alpha_k$, for $k\in[K]$, given by equations~(\ref{epsare1})
and~(\ref{alphakarealpha}), respectively, where we choose $\alpha=1$ in equation~(\ref{alphakarealpha}).
We simulate $N=400$ particles with the response function $G(s)$ and the interaction kernel $W({\mathbf x})$ 
given by
\begin{equation}
G(s) := e^{-s} 
\qquad \mbox{and} \qquad 
W({\mathbf x}) := \frac{\chi_{[0,R]}(|{\mathbf x}|)}{\pi \, R^2}\,,
\label{choiceGandW2d}
\end{equation}
where we choose the sampling radius $R=1/20=0.05$. We again discretize \eqref{SDE1} using the Euler-Maruyama 
scheme with timestep $\Delta t=10^{-3}$. We initialize the simulation by generating initial positions 
${\mathbf x}_i\in\Omega$, for $i\in[N]$, using a uniform distribution in $\Omega$, while all internal
variables are initially equal to zero, {i.e.}, ${\mathbf y}_i^k(0) = {\mathbf 0}$ for all $i\in[N]$ 
and $k\in[K]$. We again calculate $100$ stochastic realizations,
now for $K\in\{1,2,\dots,8\}$, each time performing $10^7$ timesteps.
We record the particle locations at the last timestep ($t=10^4$)
and use these for evaluating clustering properties of the model.
For identification of the clusters we again use the {\tt MATLAB} implementation
of the DBSCAN method~\cite{DBSCAN} with its parameter chosen as $\mathtt{epsilon}=0.05$ 
(same as the sampling radius) and $\mathtt{minPts}=12$.

\begin{figure}
\begin{center}
\includegraphics[width=.33\textwidth]{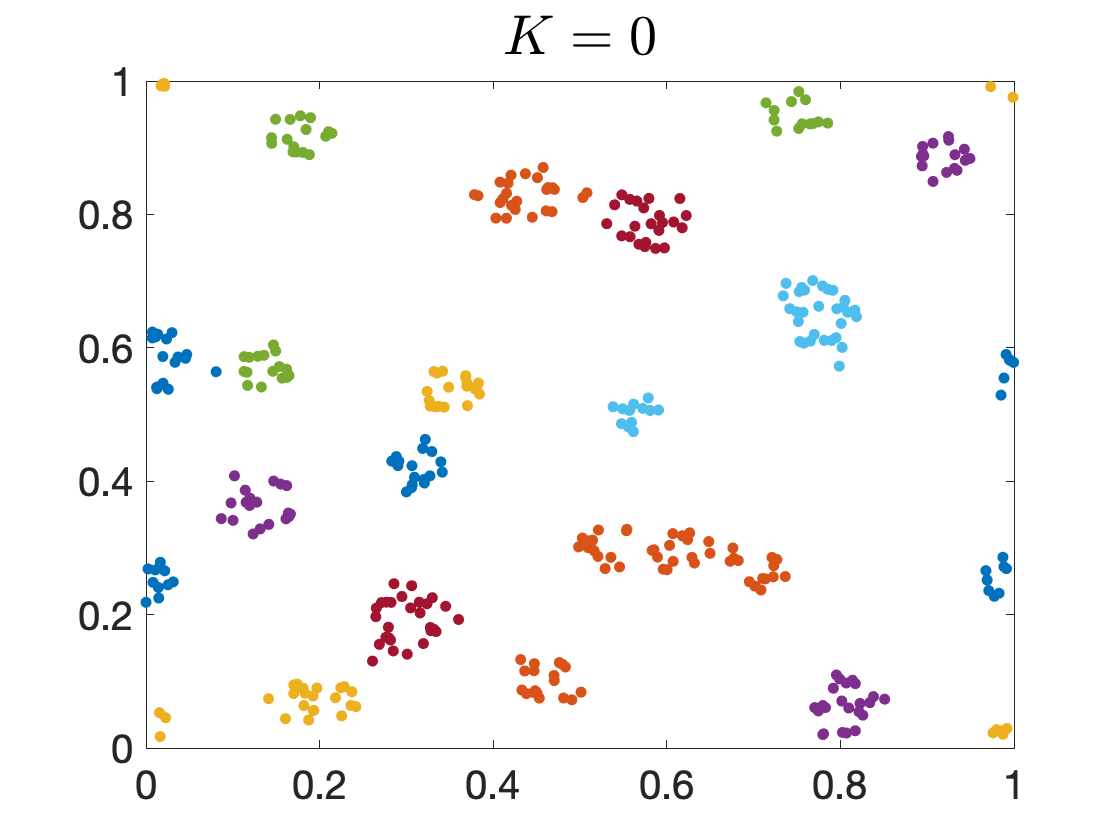}\includegraphics[width=.33\textwidth]{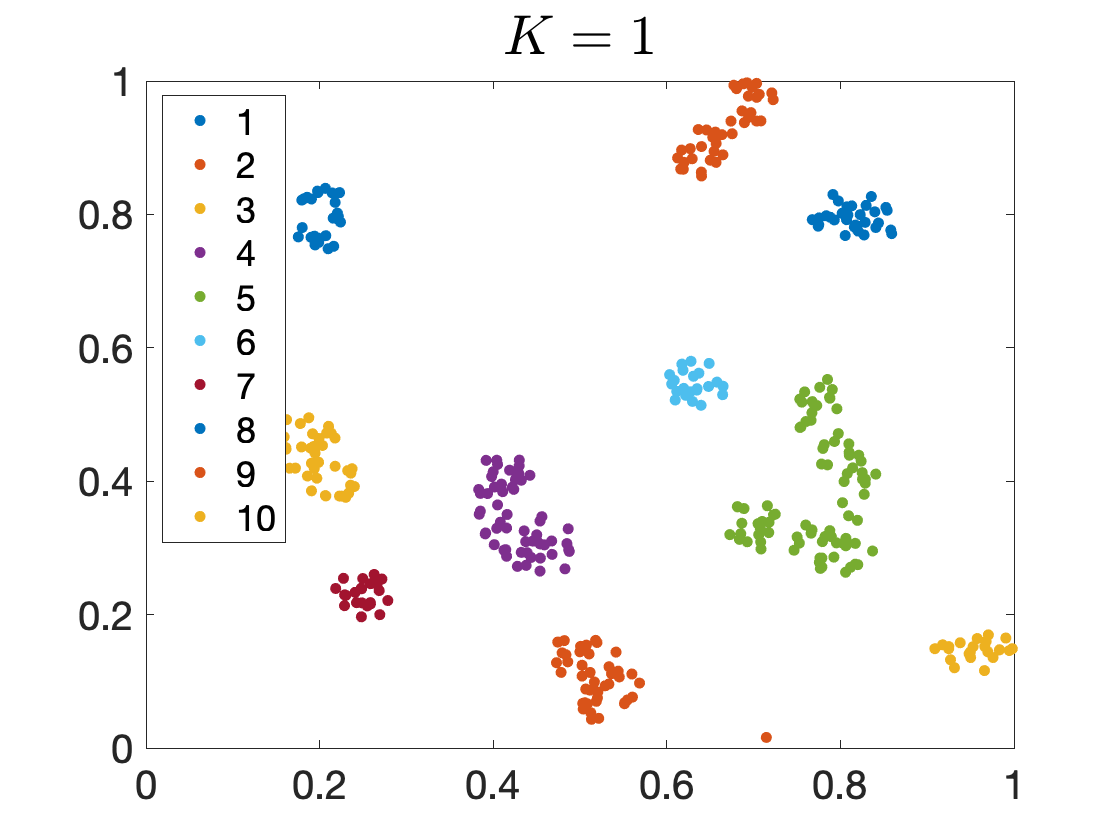}\includegraphics[width=.33\textwidth]{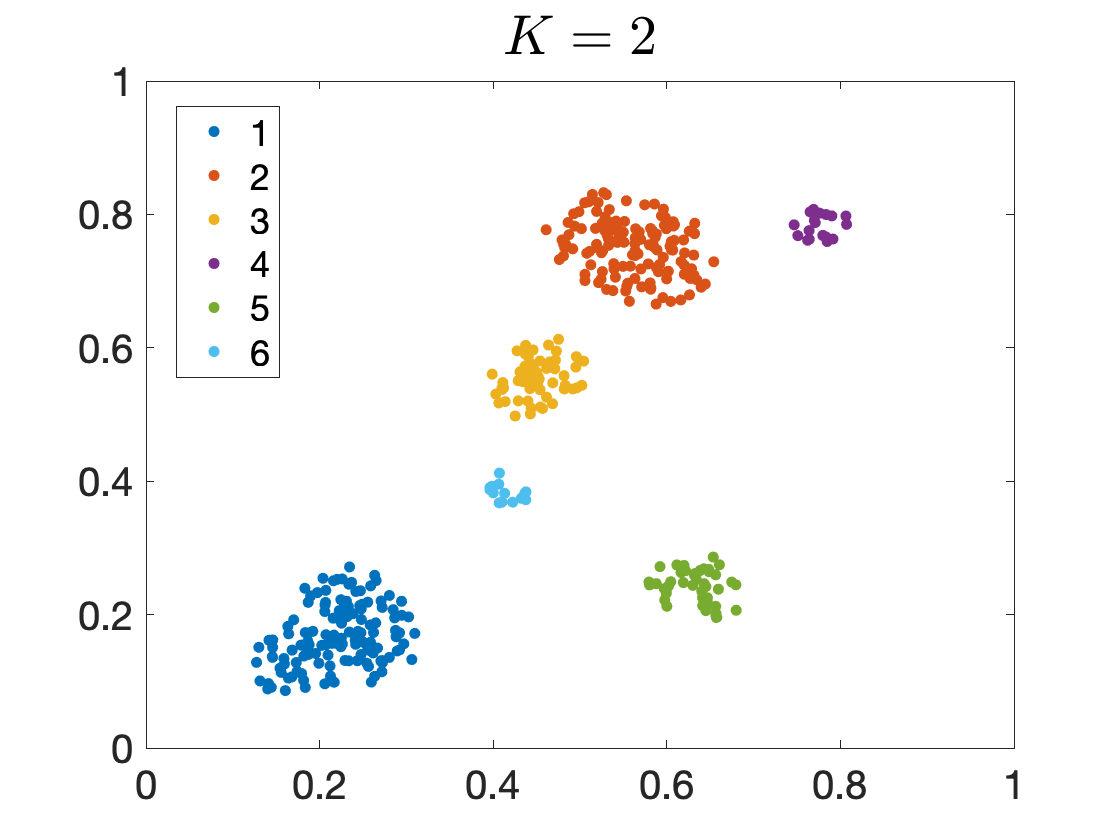} \\ 
\includegraphics[width=.33\textwidth]{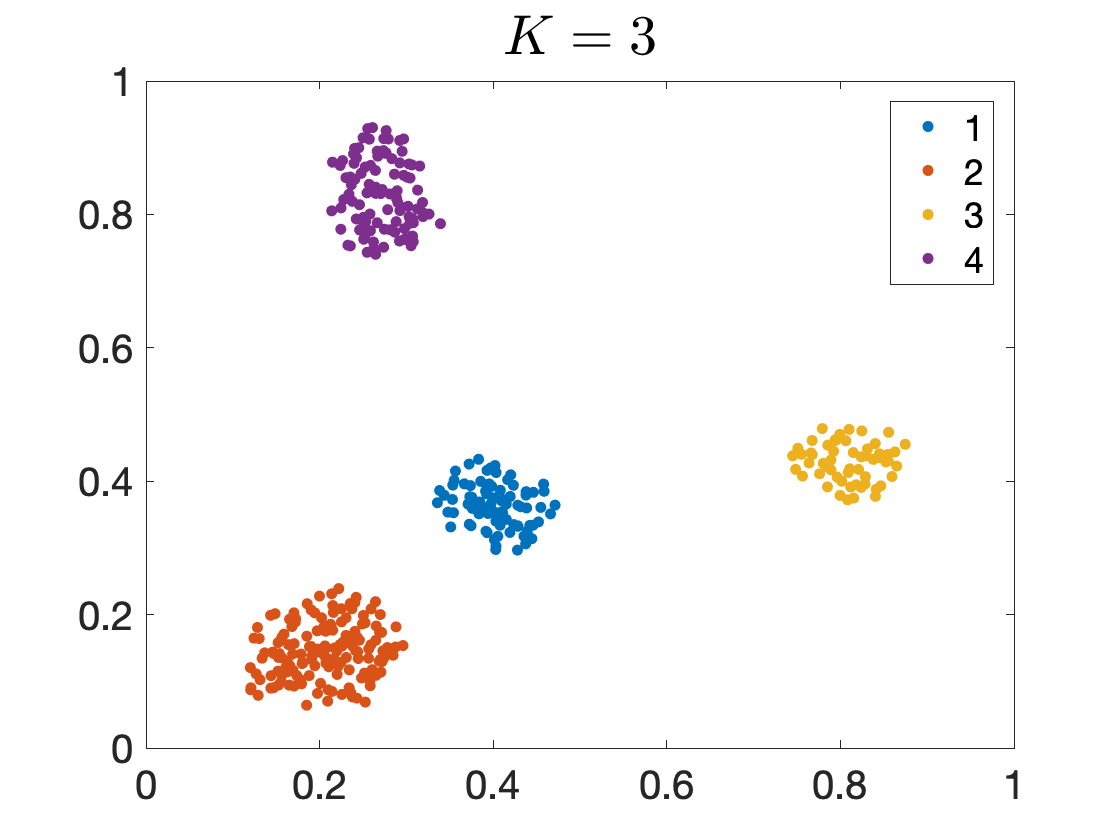}\includegraphics[width=.33\textwidth]{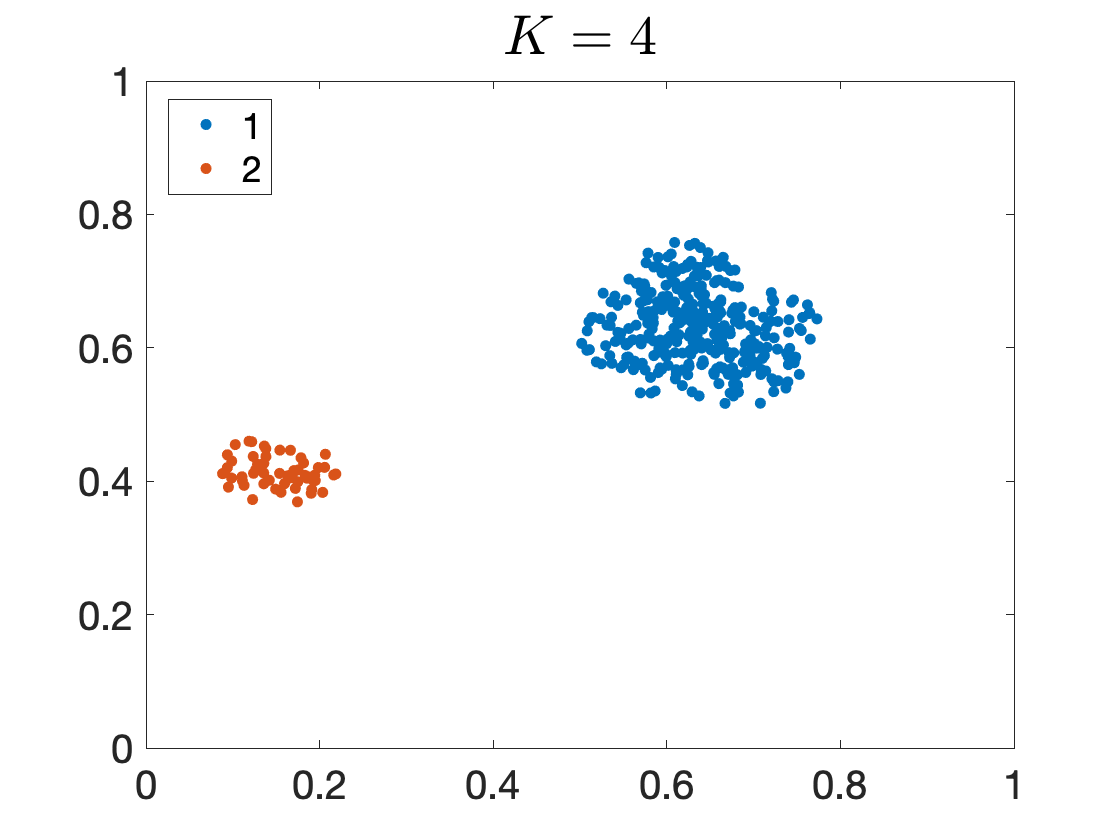}\includegraphics[width=.33\textwidth]{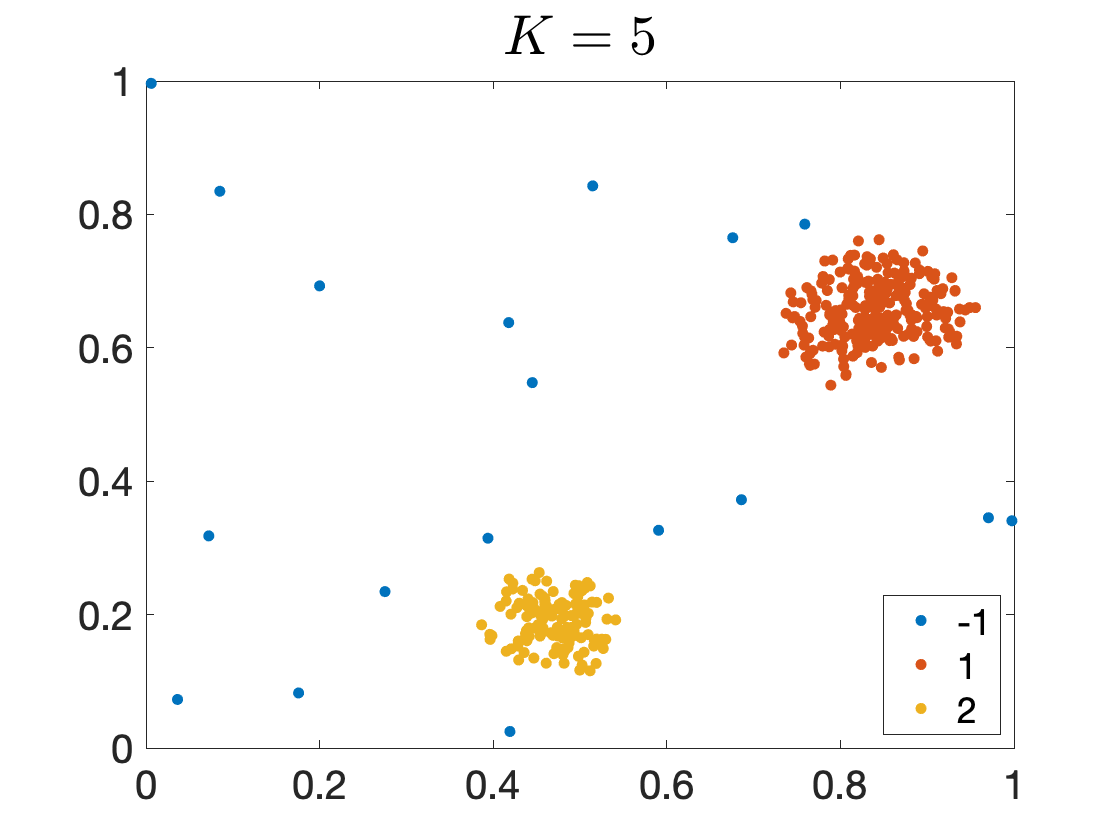} \\
\includegraphics[width=.33\textwidth]{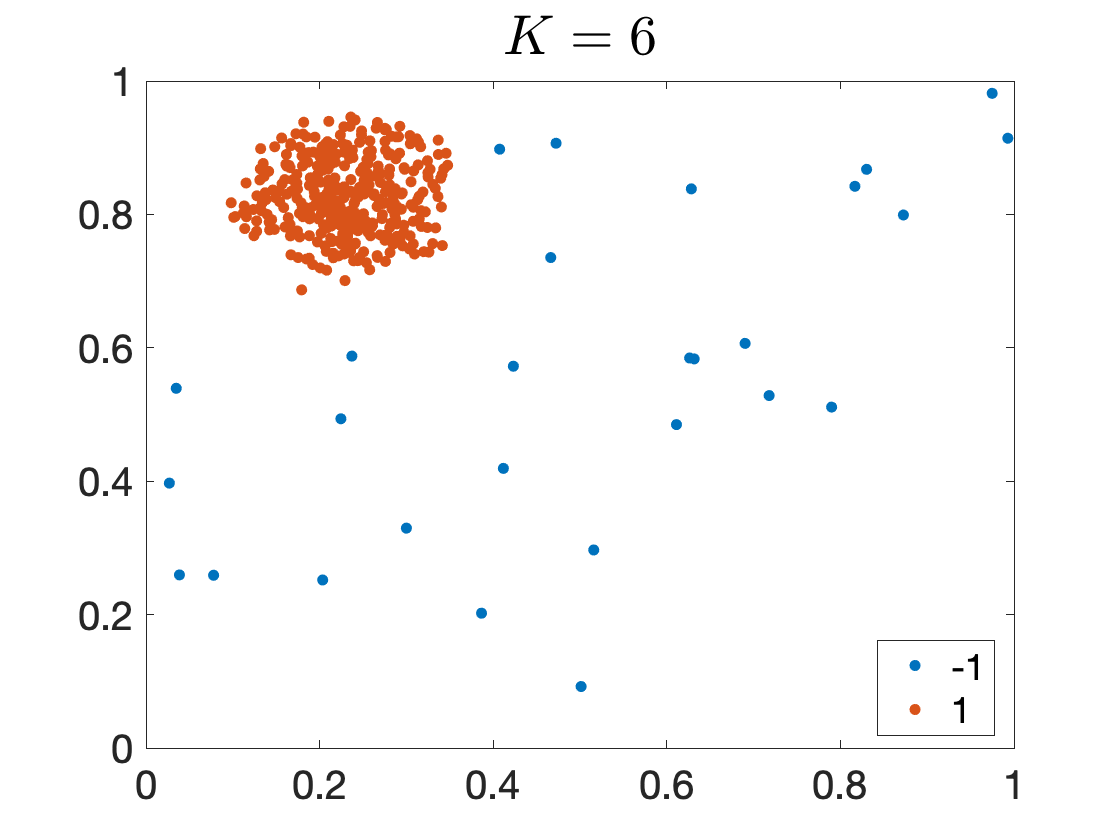}\includegraphics[width=.33\textwidth]{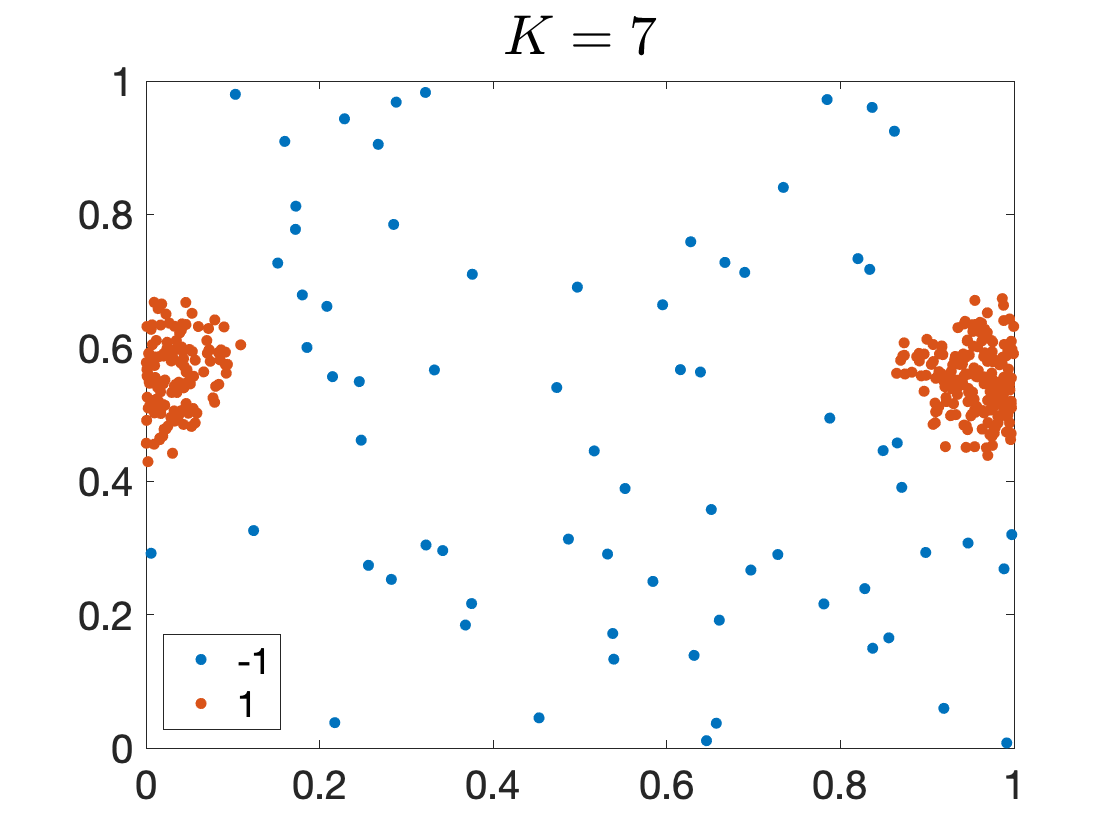}\includegraphics[width=.33\textwidth]{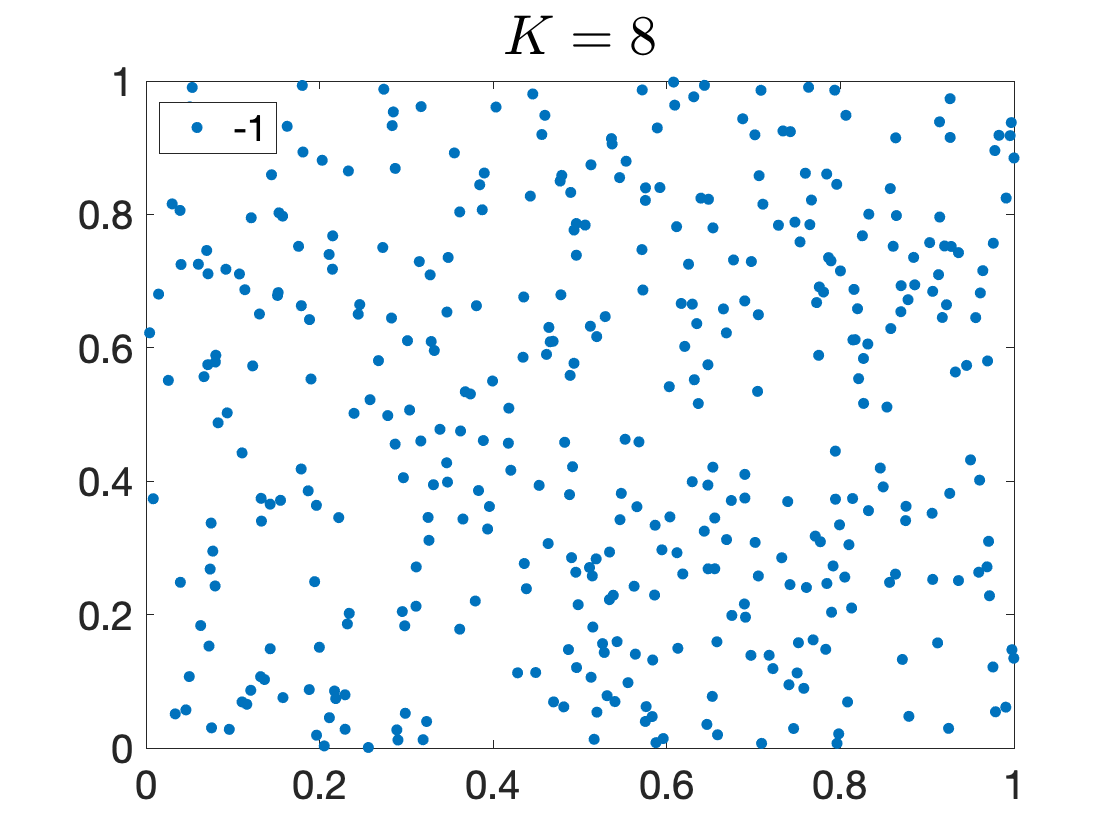}
\end{center}
\caption{
{\it 
Simulations of the individual-based model~$\eqref{rho_i}$ and~$\eqref{SDE1}$ for
$K\in\{1,\dots,8\}$, and the model~\eqref{rho_i}--\eqref{model1} without memory ($K=0)$.
We used $N=400$ agents moving in the domain $\Omega$ given by~$(\ref{domainomega})$ with $d=2$.
Functions $G(s)$ and $W({\mathbf x})$ are given by~$(\ref{choiceGandW})$ with $R=0.05$.
For each value of $K\in\{0,2,\dots,8\}$, the plots capture the particle positions at the final
time $t=10^4$. The clusters, differentiated by colour, are identified using
the {\rm DBSCAN} method with parameters $\mathtt{epsilon}=0.05$ and $\mathtt{minPts}=12$.
The points indexed with $-1$ are outliers, i.e., not belonging to any cluster.
}}
\label{figure4}
\end{figure}

Examples of the results, recorded at the final timestep of the simulations with $K\in\{1,2,\dots,8\}$,
are plotted in Figure~\ref{figure4}. The clusters are indexed by positive natural numbers, while agents 
indexed by $-1$ are classified as outliers, {i.e.}, not belonging to any cluster.
We observe the tendency to produce a smaller number of clusters with increasing $K$.
The size of clusters appears to increase until $K=4$. For larger values of $K$
the number of outliers increases. For $K=8$ no clusters are formed
in the presented stochastic realization in Figure~\ref{figure4}, 
{i.e.}, all particles have been classified as outliers.

\begin{figure}
(a) \hskip 7.5cm (b) \hfill\break
\includegraphics[width=.49\textwidth]{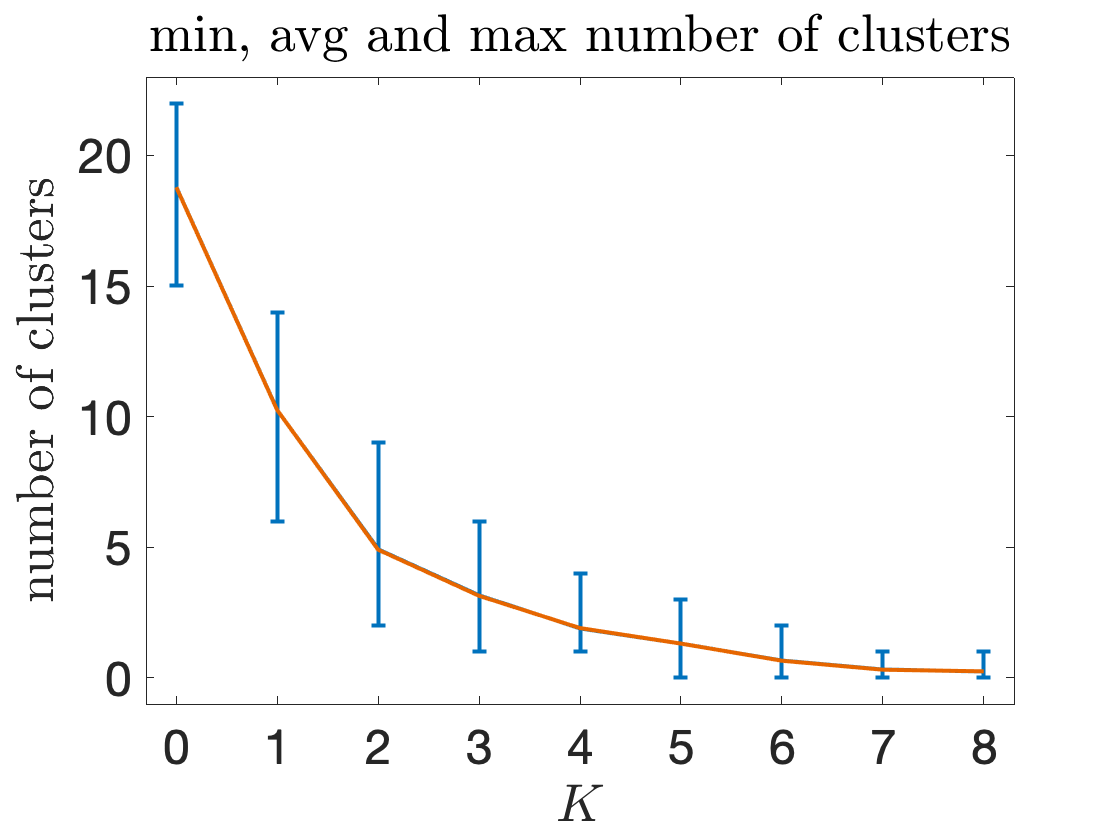}
\includegraphics[width=.49\textwidth]{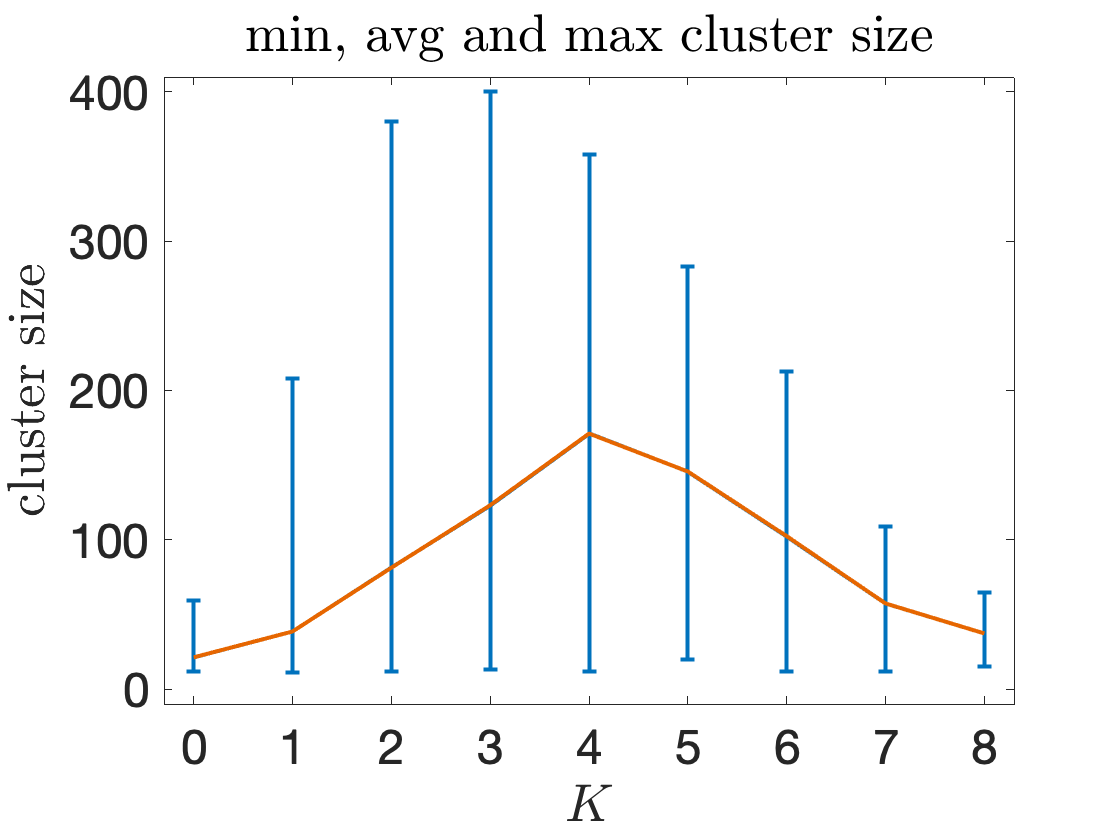} \hfill\break
(c) \hskip 7.5cm (d) \hfill\break
\includegraphics[width=.49\textwidth]{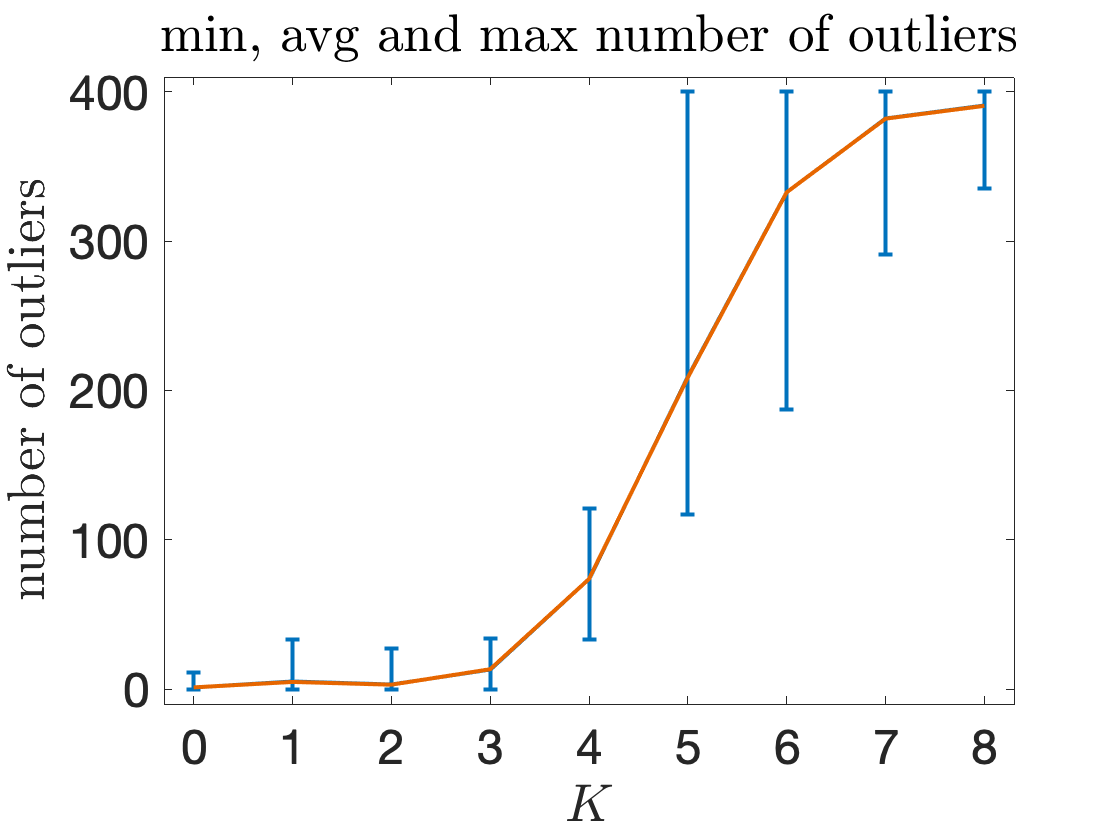}
\includegraphics[width=.49\textwidth]{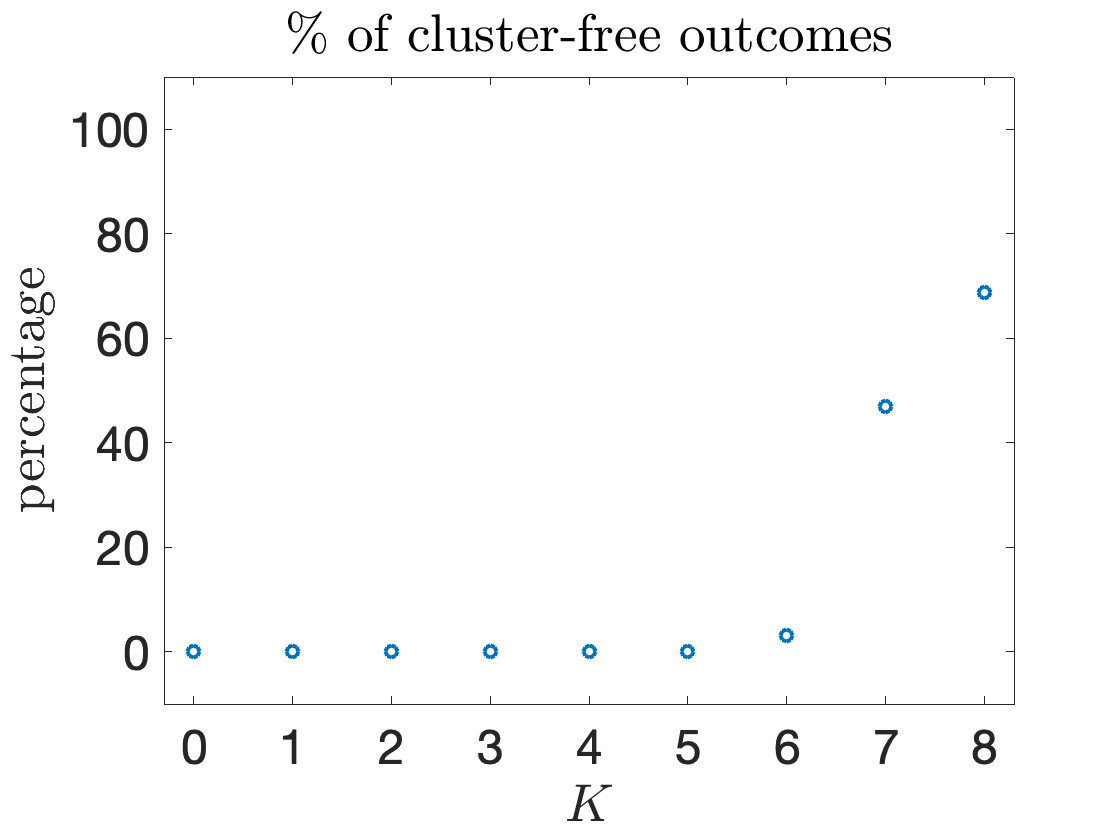}
\caption{
{\it
Statistics of the clustering behaviour over $100$ realizations of the individual-based model given 
by equations~$\eqref{rho_i}$ and~$\eqref{SDE1}$ for dimension $d=2$, $N=400$ particles, $G(s)$ and 
$W({\mathbf x})$ given by~$(\ref{choiceGandW2d})$ with $R=0.05$ and $K \in \{1, 2, \dots, 8\}$.
The case $K=0$ $($no memory$)$ refers to stochastic simulation of the system
\eqref{rho_i}--\eqref{model1}. Other parameters are the same as in Figure~$\ref{figure4}$. \hfill\break
} (a) {\it average $($orange solid line$)$, minimum and maximum $($indicated by blue error bars$)$ number of clusters identified in the 
final timestep of the simulations at time $t=10^3$, \hfill\break
} (b) {\it average $($orange solid line$)$, minimum and maximum $($blue bars$)$ cluster sizes, \hfill\break
} (c) {\it average $($orange solid line$)$, minimum and maximum $($blue bars$)$ number of outliers, {i.e.},
particles that do not belong to any cluster}, \hfill\break
(d) {\it percentage of simulation outcomes $($out of the $100$ runs$)$
that did not produce any clusters.}
}
\label{figure5}
\end{figure}

In Figure~\ref{figure5}, we provide statistics of the clustering behaviour
over the sample of $100$ realizations of the individual-based model given 
by equations~\eqref{rho_i} and~\eqref{SDE1} for $d=2$, performed for each $K \in \{1, 2, \dots, 8\}$.
Again, we also simulate the memoryless system \eqref{rho_i}--\eqref{model1}
and refer to the results by $K=0$. The statistical quantities plotted in the four panels of 
Figure~\ref{figure5} are the same as in the one-dimensional case in Figure~\ref{figure3}.
In particular, Figure~\ref{figure5}(a) shows the average, minimum and maximum number of clusters formed
obtained over $100$ stochastic realizations. We again observe the tendency to formation of fewer clusters
for higher values of $K$. For $K\in\{1,2,3, 4\}$ the size of the clusters increases
with $K$, as we observe in Figure~\ref{figure5}(b), and decreases for $K>4$.
Figure~\ref{figure5}(c) shows the average (orange solid line)
and minimum/maximum (blue error bars) number of outliers, {i.e.}, particles that do 
not belong to any cluster. We can again observe relatively low number of outliers
for $K\leq 4$ and high number for $K\geq 5$. Finally, Figure~\ref{figure5}(d) shows 
the percentage of the outcomes out of the $100$ simulations that have not produced any clusters. 
We observe that for $K\in\{1,2,3,4,5\}$ the percentage is zero (clusters have always been formed).
For larger values of $K$ the percentage increases steeply.
For $K=8$ no clusters have been formed in $75$ out of the $100$ simulations.
Based on the statistics gathered in Figure~\ref{figure5}, we may identify
three regimes:

\vskip 2mm

{\parindent -6mm \leftskip 6mm \rightskip 6mm
(a) {\it Short memory regime} ($K=1$), where a high number of clusters ($10$ on average)
is formed, typically consisting of $40$ agents. There are no or almost no outliers.

\vskip 2mm

(b) {\it Moderate memory regime} ($K=4$), where typically a few larger clusters form, containing
a significant proportion of agents except for the outliers.

\vskip 2mm

(c) {\it Long memory regime} ($K=8$) where most of the simulations do not produce
any clusters, and even if a cluster is formed, the majority of agents exists as outliers.
\par}

\vskip 2mm

\begin{figure}
(a) \hskip 7.5cm (b) \hfill\break
\includegraphics[width=.49\textwidth]{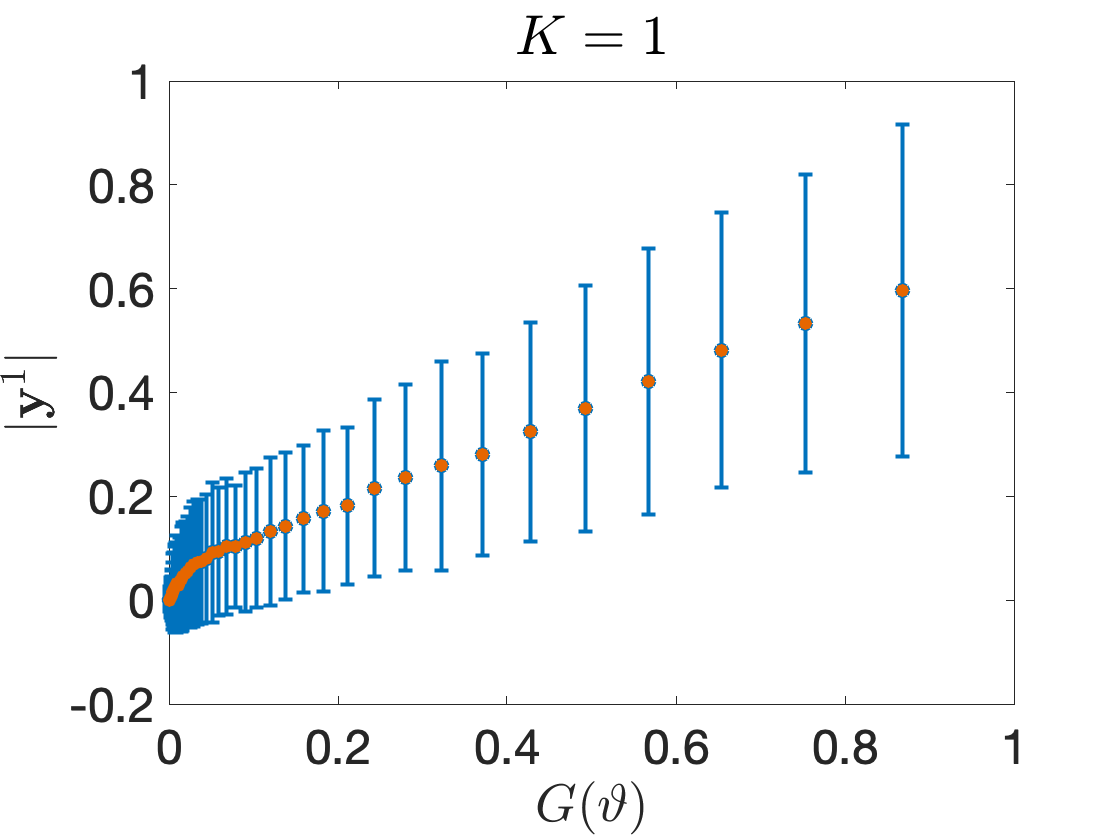}\includegraphics[width=.49\textwidth]{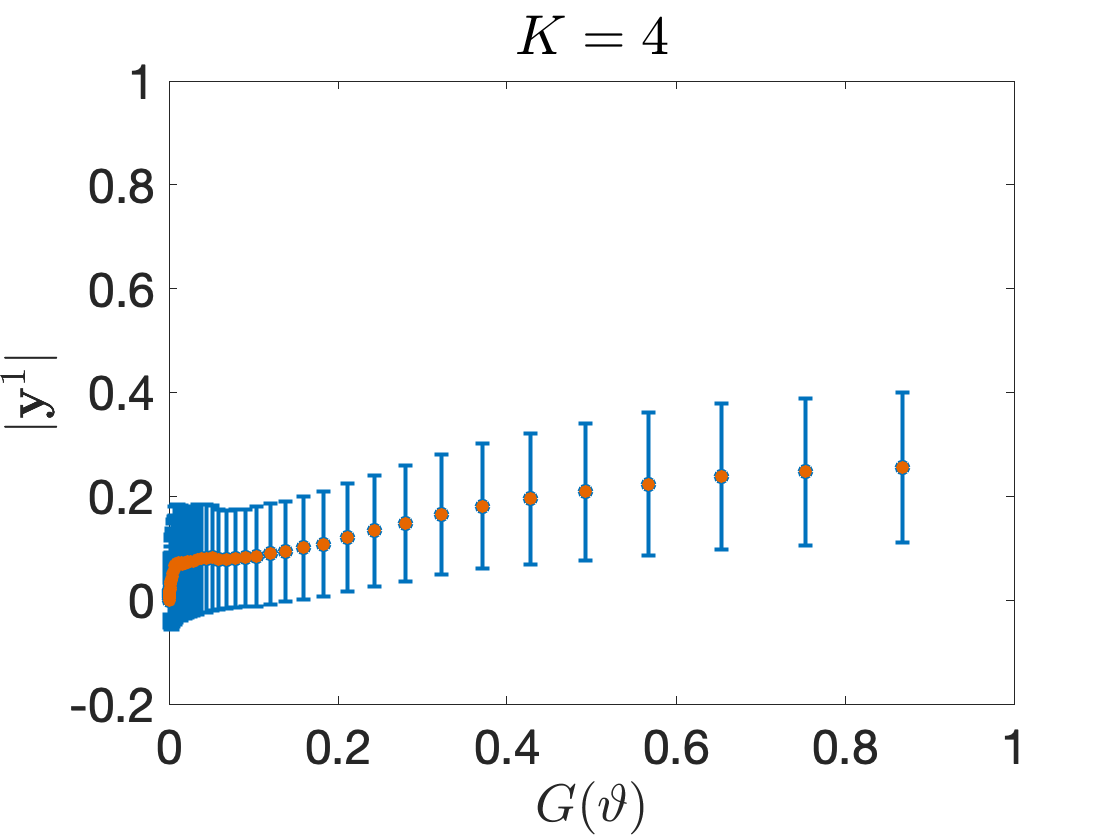}\hfill\break
(c) \hskip 7.5cm (d) \hfill\break
\includegraphics[width=.49\textwidth]{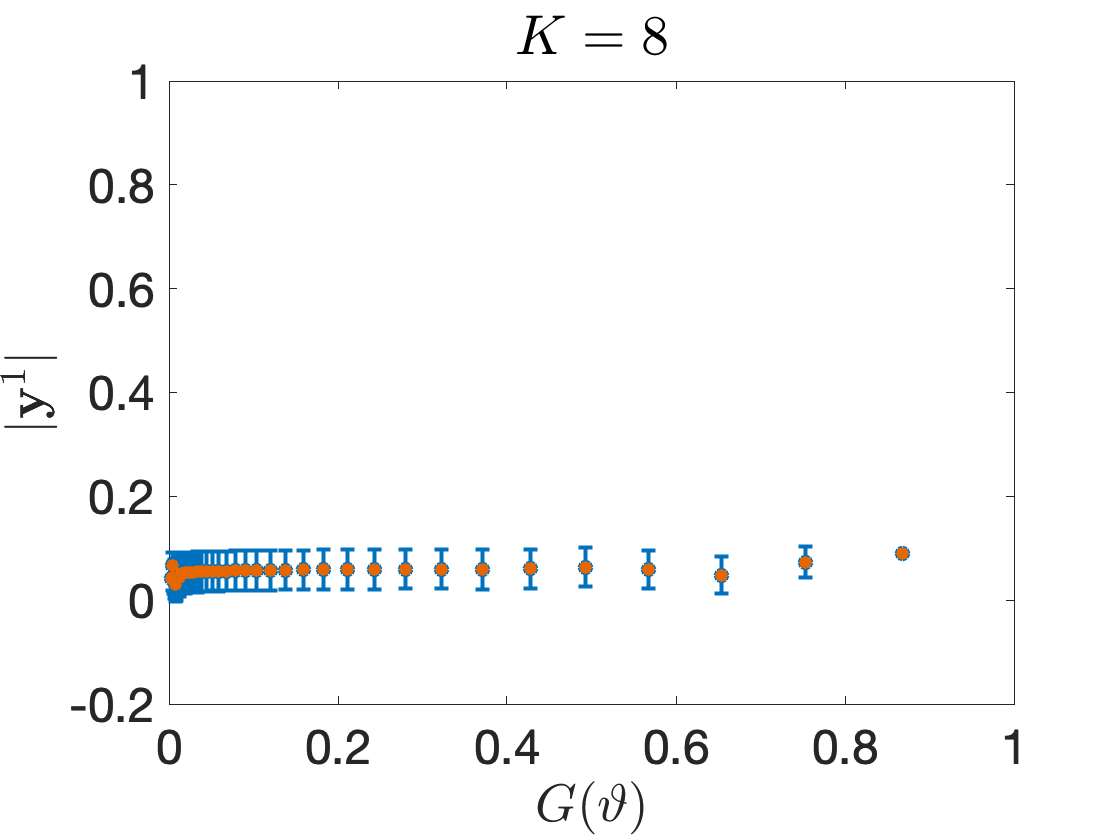}\includegraphics[width=.49\textwidth]{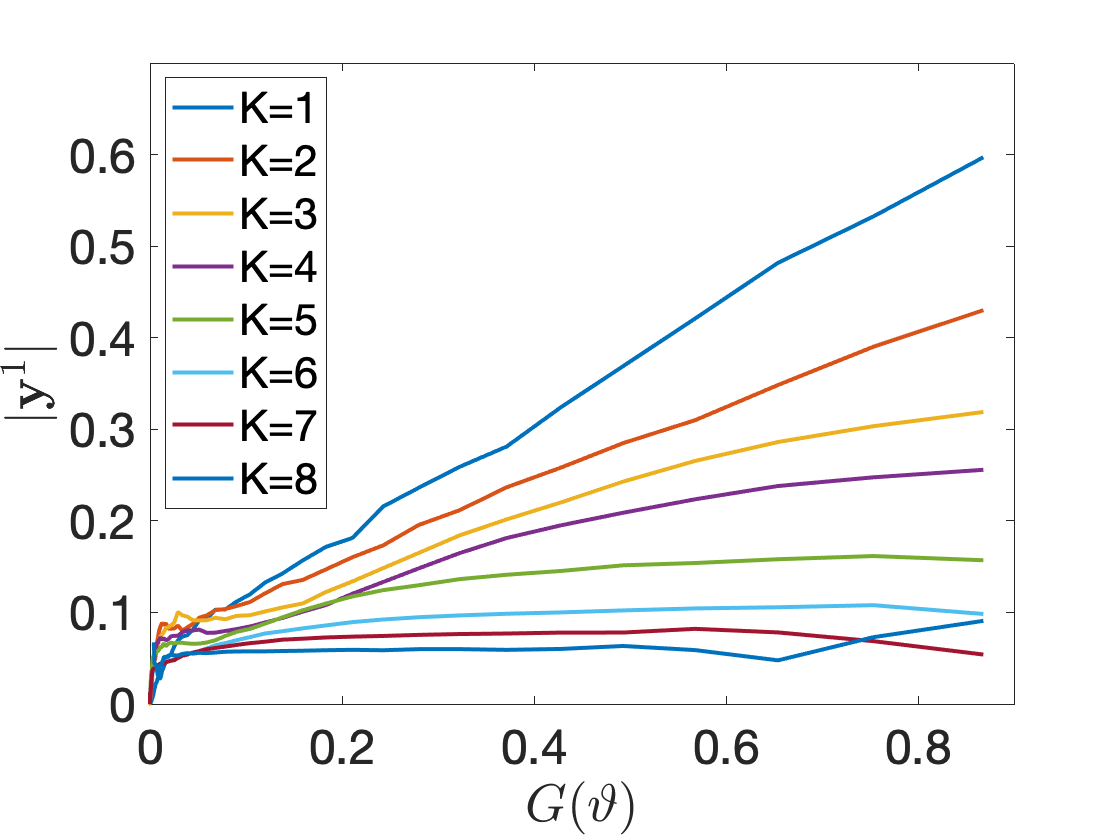}
\caption{{\it
Statistics of the magnitude of the internal variable $|{\mathbf y}_i^1|$ versus the response function value $G(\vartheta_i)$,
aggregated over all agents and timesteps in $100$ simulations of the individual-based model given by equations~$\eqref{rho_i}$ 
and~$\eqref{SDE1}$ in $d=2$ dimensions with $N=400$ particles,
for:} \hfill\break
(a) $K=1$ ({\it short memory regime}), \hfill\break
(b) $K=4$ ({\it moderate memory regime}), \hfill\break
(c) $K=8$ ({\it long memory regime}). \hfill\break 
{\it The orange circles represent the average value of $|{\mathbf y}_i^1|$ corresponding to each $G(\vartheta_i)$, while the error bars indicate 
the standard deviations. The panel} (d) {\it presents the average values of $|{\mathbf y}_i^1|$ corresponding to each $G(\vartheta_i)$ 
for all values $K \in\{1,2,\dots,8\}$.
}}
\label{figure6}
\end{figure}

\noindent
To gain further understanding of this behaviour, we have collected the statistics
of the values of the internal variables ${\mathbf y}^1_i$, for $i\in[N]$,
recorded during the temporal evolution of the individual-based model
in the considered $100$ stochastic realizations in Figure~\ref{figure6}.
More precisely, the forward Euler discretization of the first equation 
of the SDE system~\eqref{SDE1} reads
$$
{\mathbf x}_i(t+\Delta t) \,=\, {\mathbf x}_i(t) \,+\, {\mathbf y}_i^1(t) \, \Delta t \qquad\quad \mbox{for} \quad i\in[N]\,,
$$
{i.e.}, in each timestep of the discrete simulation the particle
locations in the physical space are updated by ${\mathbf y}_i^1(t)\, \Delta t$.
In Figure~\ref{figure6}, we plot the statistics of the magnitudes $|{\mathbf y}_i^1(t)|$
in dependence of the values of $G(\vartheta_i(t))$ with
the perceived densities $\vartheta_i(t)$ given by equation~\eqref{rho_i}.
In particular, for each fixed $K\in[8]$, in each timestep of each of the $100$
stochastic realizatons, we record the pair $(G(\vartheta_i(t)), |{\mathbf y}_i^1(t)|)$
for all $i\in[N]$. We then plot the mean and standard deviation of
the values of $|{\mathbf y}_i^1|$ against the values of $G(\vartheta_i)$,
separately for each memory regime in panels (a), (b) and (c) of Figure~\ref{figure6}.
Figure~\ref{figure6}(a) corresponds to the short memory regime, Figure~\ref{figure6}(b)
presents the moderate memory regime and Figure~\ref{figure6}(c) corresponds
to the long memory regime. We observe that $|{\mathbf y}_i^1|$ increases 
with $G(\vartheta_i)$, which corresponds to the fundamental modeling assumption
that the particles' mobility decreases with perceived density
(recall that the response function $G$ is in general a decreasing function;
we use $G(s) = e^{-s}$ in equation~(\ref{choiceGandW2d})).
However, this effect weakens as $K$ increases, and is nearly absent at $K=8$,
see Figure~\ref{figure6}(c). We therefore conclude that the presence of memory
systematically inhibits the particles' responsiveness to environmental stimuli,
specifically the perceived density of their neighbours.

\section{Discussion}

\label{sec7}

In this paper, we have shown that short-term memory enhances and long-term memory inhibits spontaneous particle 
aggregation. Our investigation has been based on the (first-order) spontaneous aggregation model
without memory which has been previously investigated in the literature~\cite{BHW:2012:PhysD}. Its main properties 
are summarized in Section~\ref{sec2}. The memory has been added into this model in Section~\ref{sec3} by introducing 
a chain of $K$ internal variables that allow the agents to `remember' the densities they encountered in the past. 
If $K=1$, then our individual-based model is equivalent to what is called  `the second-order model' in reference~\cite{BHW:2012:PhysD}, 
where the internal variable represents the agent's velocity. Considering $K>1$, the model introduces additional
layers of memory described by $K$ internal variables. Additional internal variables can be used to better
fit the properties of relatively complex (high-dimensional) models of interacting particle systems, while keeping the number of degrees
of the system (internal variables) relatively small~\cite{Erban:2016:CAM,Erban:2020:CGM}. For simplicity, our transfer of 
information between layers of memory is linear, but nonlinear functions can also be introduced to better fit the properties of 
some systems~\cite{Erban:2020:CGM, Gurney:2018:INN, Erban:2025:arxiv}.

In Sections~\ref{sec5} and~\ref{sec6}, we have reported the results of 
systematic stochastic simulations of the individual-based model in the spatially one-dimensional 
and two-dimensional settings, respectively. We have shown that memory inhibits the particles' 
responsiveness to environmental stimuli. Our results show that the introduction of memory leads, in general, 
to the formation of a smaller number of larger clusters. This trend is observed until $K=3$ in the one-dimensional 
setting and until $K=4$ in two spatial dimensions. When the number of layers $K$ is increased further, 
{i.e.}, as the memory becomes `longer', its effect starts to be disruptive. This is manifested by the increasing 
proportion of `outliers', which are the particles that are not part of any cluster. Also, the percentage of 
stochastic simulations where no clusters are formed during the observation time increases with
increasing $K$. We therefore conclude that short-term or medium-term memory has a coarsening
effect on spontaneous aggregation, while long-term memory disrupts it.

The presence of memory impacts the clustering properties of the model - in particular, the number and size 
of clusters being formed, or even the very ability of the model to produce clusters, starting from initially 
randomly distributed particles in the physical space. In one spatial dimension, we have observed a sharp 
transition between $K\leq 4$, when clusters are always formed, and $K\geq 5$, when clusters are almost never 
formed during the simulation period. A similar, although less sharp transition, takes
place in two-dimensional models between $K\leq 6$ and $K\geq 7$.

We have performed systematic numerical simulations to evaluate the impact of the number of internal 
variables (or `memory layers') $K$ on the dynamics of the many particle system, because (exact) analytical
mathematical results to study these properties are not available and the development of the corresponding
mathematical theories includes a number of open questions. Some progress in this direction can be made as shown
in Section~\ref{sec:FP}, where we have derived the formal macroscopic limit of the system
as the number of agents tends to infinity, under the usual molecular chaos assumption~\cite{BHW:2012:PhysD,Erban:2012:ICB}.
The limit is described by a Fokker-Planck equation, and we have characterized its steady states to gain an insight 
into the patterns (clusters) formed by the system. 

Many organisms benefit from learning to adapt their behaviour to the distribution of resources and
their learning depends on the memory capacity of each individual~\cite{Falcon:2023:LMS}. In our 
investigation, the number of memory layers, $K$, has been fixed, but it could be made variable
for each agent. If the evolutionary objective is to develop large clusters in two-dimensional setting,
then the choice $K=4$ can be considered `optimal'. Our memory model with $K$ layers could also be
further extended by considering nonlinear activation functions to model a more complicated
neural network~\cite{Gurney:2018:INN} and the applicability of the model could be further enhanced 
by introducing different modes of learning during the individual-based simulations~\cite{Lewis:2021:LAM}.

\rev{In our computational simulations, we have used exponentially decaying function $G(s) = e^{-s}$ in 
equations~(\ref{choiceGandW}) and~(\ref{choiceGandW2d}) in order for clusters to be formed on effectively
`intermediate' timescales of the order $10^3$. Since all of the simulated particles are subject to the Brownian motion
(including the particles in the clusters), the clusters are also able to slowly move (diffuse) on much longer 
timescales. However, as we are interested in biological applications, such long time behaviour of the model
can be irrelevant, because other processes (e.g birth or death of individuals) could also happen on long
timescales. We have therefore run the simulation until clusters emerge (starting from the initial random distribution 
of particles) and become `quasi-stationary' in the sense that they do not move (on this timescale) and are 
in dynamic equilibrium with the outliers. For a study of the clustering behaviour on long timescales, see, e.g.,~\cite{Carrillo:2019:AP}.}




\section*{Declarations}

\subsection*{Authorship and Contributorship}
R.E. and J.H. contributed equally to the conception and writing 
of this manuscript.
Both authors have read and approved the final manuscript.
In addition, the following contributions occurred:
J.H. developed the computer code, performed the numerical simulations, and prepared all figures.

\subsection*{Conflict of interests}
The authors declare no conflict of interests.

\subsection*{Data \& Code Availability}
In compliance with EPSRC's open access initiative, the code used in this paper is available from
\href{https://doi.org/10.5281/zenodo.18135482}{\tt doi.org/10.5281/zenodo.18135482}.

\subsection*{Ethics}
This study did not involve human participants, animals, or sensitive data requiring ethical approval or consent to participate.

\subsection*{Funding}
This work was supported by the Engineering and Physical Sciences Research Council, grant number EP/V047469/1, awarded to Radek Erban.

\section*{Acknowledgements}
The authors thank to Martin Burger (DESY and Universit\"{a}t Hamburg) for stimulating discussions that contributed to the preparation of this paper.


\bibliography{bibmemory}

\end{document}